\documentclass[12pt,usenames,dvipsnames]{amsart}

\DeclareFontFamily{U}{matha}{\hyphenchar\font45}
\DeclareFontShape{U}{matha}{m}{n}{
	<5> <6> <7> <8> <9> <10> gen * matha
	<10.95> matha10 <12> <14.4> <17.28> <20.74> <24.88> matha12
}{}
\DeclareSymbolFont{matha}{U}{matha}{m}{n}
\DeclareMathSymbol{\odiv}         {2}{matha}{"63}

\usepackage[margin=0.8in]{geometry}

\usepackage{color}
\usepackage{fancyhdr}
\usepackage[utf8]{inputenc}
\usepackage[T1]{fontenc}
\usepackage{lmodern}
\usepackage{amsmath,amssymb}
\usepackage[english]{babel}

\usepackage{mathtools}
\usepackage{mathrsfs}
\usepackage{enumitem}
\usepackage{amsthm}
\usepackage{amsfonts}
\usepackage{stmaryrd}
\usepackage{mathdots}

\usepackage{graphicx}
\usepackage{verbatim}
\usepackage{dsfont}
\usepackage{bbm}

\usepackage{caption}
\usepackage{subcaption}
\usepackage{tikz-cd}
\usepackage{tikz}
\usetikzlibrary{arrows, positioning, calc, shapes}

\usepackage{fancyhdr}
\usepackage{appendix}
\usepackage{csquotes}

\usepackage[numbers]{natbib}
\usepackage[english]{hyperref}

\graphicspath{{images/}}

\newcommand\R{\mathbb R}
\newcommand\Z{\mathbb Z}

\newcommand\C{\mathbb C}
\renewcommand\P{\mathbb P}

\newcommand{\cyc}[2]{\left\llbracket #1 , #2 \right\rrbracket}

\definecolor{darkred}{rgb}{0.7, 0.0, 0.0}
\definecolor{darkblue}{rgb}{0.0, 0.2, 0.6}
\definecolor{tblue}{rgb}{0, 0.53333,1}
\definecolor{tgray}{rgb}{0.9, 0.9, 0.9}

\newtheorem{theorem}{Theorem}[section]
\newtheorem{proposition}[theorem]{Proposition}

\newtheorem{lemma}[theorem]{Lemma}
\newtheorem{convention}[theorem]{Convention}
\newtheorem{conjecture}{Conjecture}[section]

\newtheorem{definition}[theorem]{Definition}

\theoremstyle{remark}
\newtheorem{remark}[theorem]{Remark}

\theoremstyle{remark}
\newtheorem{example}[theorem]{Example}

\newcommand{\xqed}[1]{
	\leavevmode\unskip\penalty9999 \hbox{}\nobreak\hfill
	\quad\hbox{\ensuremath{#1}}}


\title{Newton-Okounkov bodies obtained from certain orbits of plabic graphs}
\author{Michael Schl\"o{\ss}er}
\date{}

\begin{document}
	
	\begin{abstract}
		We investigate the plabic graphs corresponding to the quadrilateral Postnikov arrangements used by J.Scott to equip the homogeneous coordinate rings of Grassmannians with a cluster structure. More precisely we describe their orbits under the natural action of the dihedral group and show that the associated Newton-Okounkov bodies are all unimodular equivalent to Gelfand-Tsetlin polytopes.
	\end{abstract}
	
	\maketitle
	
	\section{Introduction}
	
	Newton-Okounkov bodies are defined as closed convex hulls of points and even when they are polyhedral it is in general quite hard to obtain their defining inequalities. 
	In their paper \textit{Cluster duality and mirror symmetry for Grassmannians} \cite{Rietsch_2019}, K. Rietsch and L. Williams construct Newton-Okounkov bodies for Grassmannians and use mirror symmetry to obtain defining inequalities from the tropicalization of the super-potential on an open set of the mirror Grassmannian arising from the Landau-Ginzburg model. 
	In particular they associate Newton-Okounkov bodies to certain planar graphs going by the name \textit{plabic graphs}.
	These were introduced by A. Postnikov in \cite{Postnikov2006} to study the totally nonnegative part of the Grassmannian and are defined as planar graphs drawn inside a disk, with boundary vertices labeled by 1, \ldots, n and internal vertices coloured either black or white. \\
	
	While all Newton-Okounkov bodies coming from plabic graphs turn out to be polytopes \cite{Rietsch_2019}, some families of plabic graphs lead to especially nicely behaved bodies. 
	For example the Newton-Okounkov body associated to \textit{rectangle graph} $G^\text{rec}_{k,n}$ \cite{Rietsch_2019} (see section 3) can be shown to be unimodular equivalent to certain Gelfand-Tsetlin polytope $\mathcal{GT}_{k,n}^1$ \cite{Rietsch_2019}.
	In their paper \textit{symmetries on plabic graphs and associated polytopes} \cite{FF2018}, X. Fang and G. Fourier on the other hand show that the Newton-Okounkov bodies of what they call \textit{dual rectangle graphs} (see section 3) are equivalent to FFLV-polytopes, originally introduced in \cite{Feigin2011}.
	While both Gelfand-Tsetlin and FFLV-polytopes are integral, one can in general not expect to obtain an integral polytope when starting from an arbitrary plabic graph. Indeed, in \cite{Rietsch_2019}, \cite{Rappel2017} examples of plabic graphs are given, which lead to non integral polytopes. In particular they are neither unimodular equivalent to Gelfand-Tsetlin nor FFLV-polytopes. \\
	
	In this paper we describe two additional families of plabic graphs, which we call \textit{checkboard graphs} and \textit{dual checkboard graphs} respectively (see Figure \ref{fig:figure3}). They correspond to the quadrilateral arrangements used by Scott in \cite{Scott2006} to equip the homogeneous coordinate rings of Grassmannians with a cluster structure.
	Checkboard graphs possess an especially simple internal structure. Indeed, unlike most other plabic graphs, all of their non boundary faces are all given by squares.
	
	The faces of any given plabic graph may be labelled by $k$-subsets of $\{1, \ldots, n\}$. These labels form weakly separated collections (see section 2) and uniquely determine the corresponding plabic graph (up to some local modifications) \cite{Oh2015}. Weak separation is easily seen to be preserved under the natural action of the dihedral group $D_n$ which we view as subgroup of the symmetric group $S_n$. We also prove the converse.
	
	\begin{proposition}
		Let $\rho \in S_n$, then $\rho$ preserves weak separation if and only if $\rho \in D_n$. In particular $D_n$ acts on the set of maximal weakly separated collections.
	\end{proposition}
	
	We thus obtain an action of $D_n$ on plabic graphs by identifying them with their labels. 
	This action is compatible with the various local transformations called \textit{moves} on plabic graphs (see section 4). 
	Naturally we explore the orbits of the aforementioned families of graphs. If $n=4$ there is only one orbit containing $2$ graphs. Assuming $n>4$, the main takeaways may be summarized as follows.

	\begin{proposition}
		Neither the checkboard graph nor its dual lies in the orbit of a rectangle graph.
		The dual rectangle graph $\widehat G^\text{rec}_{k,n}$ is a reflection of the rectangle graph $G^\text{rec}_{k,n}$ and its stabilizer under $D_n$ is trivial.
		The dual checkboard graph $\widehat G^\text{ch}_{k,n}$ is a rotation of the checkboard graph $G^\text{ch}_{k,n}$ if either $k$ or $n-k$ is even. 
		Otherwise these graphs lie in distinct orbits. In any case the union of their orbits contains $n$ graphs i.e.
		$|D_n G^\text{ch}_{k,n} \cup D_n \widehat G^\text{ch}_{k,n}| = n$.
	\end{proposition}
	
	Finally we calculate the Newton-Okounkov bodies for orbits of (dual) checkboard graphs. 
	
	\begin{theorem}
		Let $G$ be any graph in the $D_n$-orbit of the checkboard graph $G^\text{ch}_{k,n}$ or the dual checkboard graph $\widehat G^\text{ch}_{k,n}$. Then its associated  Newton-Okounkov body $\Delta_G$ is unimodular equivalent to the Gelfand-Tsetlin polytope $\mathcal{GT}_{k,n}^1$.
	\end{theorem}
	
	A major step in proving this result is deriving closed formulas expressing the super-potential in terms of Plücker coordinates associated to $G \in D_n G^\text{ch}_{k,n} \cup D_n \widehat G^\text{ch}_{k,n}$. Essentially this boils down to identifying sequences of moves transforming $G$ into graphs containing certain target labels. 
	This approach more elementary than the proof of the analogous result for rectangle graphs (see \cite{Rietsch_2019} Proposition 10.5, \cite{Marsh2020} Proposition 6.10). \\
	
	As a consequence of Proposition 1.2 and the results of \cite{Rietsch_2019} as well as \cite{FF2018}, we see that acting with a reflection on a plabic graph may result in a Newton-Okounkov body which is not unimodular equivalent to the one associated to the original graph. However we make the following conjecture.
	
	\begin{conjecture}
		All rotations of a plabic graphs give the same Newton-Okounkov body (up to unimodular equivalence).
	\end{conjecture}
	
	Dealing with plabic graphs and their related combinatorics by hand can be quite tedious. This motivated the development of a Julia package implementing
	much of the combinatorics related to plabic graphs. Among other things this includes tools for interactive visualization and plotting. Moreover various searching algorithms are provided and the calculation of Newton-Okounkov bodies associated to plabic graphs is automated. The package together with a comprehensive documentation can be found at \cite{Schloesser2024}. \\

	This paper is split into 10 sections. After the introduction we recall plabic graphs and introduce the relevant families of plabic graphs in sections 2 and 3. The action of the dihedral group on plabic graphs, in particular orbits and stabilizers of rectangle and checkboard graphs are investigated in section 4. Sections 5 to 7 are devoted to recalling quiver mutations, Grassmannians and their $\mathcal{A}$-cluster structure. Then in section 8 we briefly recall how to obtain Newton-Okounkov bodies from tropicalizing the superpotential, before proving our main results in section 9 and 10.
	
	\section*{Acknowledgements}
	This paper is based on the author's master thesis and he would like to thank his supervisors Ghislain Fourier and Xin Fang for many helpful discussions and input during the development of this paper.
	The author was supported by the Deutsche Forschungsgemeinschaft
	(DFG, German Research Foundation), Project-ID 286237555, TRR 195.
	
	\section{Plabic graphs}
	We recall the definition of plabic graphs and some of their properties following \cite{Postnikov2006} and \cite{Oh2015} (up to some minor technicalities). See also \cite{Oh_2020}. \\
	
	Let us quickly set up some notation. 
	Throughout this paper, if not stated otherwise, $k, n$ will always denote integers such that $2 \leq k \leq n-2$.
	Moreover for any integer $m \geq 1$ we use the notation $[m]:= \{1, 2, \ldots, m\}$ and furthermore set $[0] = \emptyset$. 
	The set of subsets of cardinality $k$ inside $[n]$ is denoted by $\binom{[n]}{k}$. We commonly identify $[n]$ with $\Z / n \Z$, thus giving $[n]$ a natural cyclic order. Finally we denote cyclic intervals by
	$$ \cyc{a}{b} = \begin{cases}
		\{a, a+1, ..., b\} & \text{ if } a \leq b \\ \{a, a+1 , ..., n\} \cup \{ 1, 2, ... ,b\} & \text{ if } a > b.
	\end{cases} .$$
	
	\begin{definition}
		A \textbf{plabic graph} is a finite connected plane graph $G$ whose interior is bounded by a vertex disjoint cycle containing $n$ \textbf{boundary vertices} $b_1, \ldots, b_n$. Here the ordering is chosen to be clockwise. We also use the notation $i$ for $b_i$. 
		The vertices (edges) in the interior of $G$ are also called \textbf{inner} vertices (edges). Any boundary vertex should be incident to exactly one inner edge and any inner vertex is coloured either black or white.
		We consider plabic graphs up to orientation preserving homeomorphisms.
	\end{definition}
	\begin{figure}[!htb]
		\centering
		\begin{tikzpicture}
			\node (background) {\includegraphics[height=6.5cm]{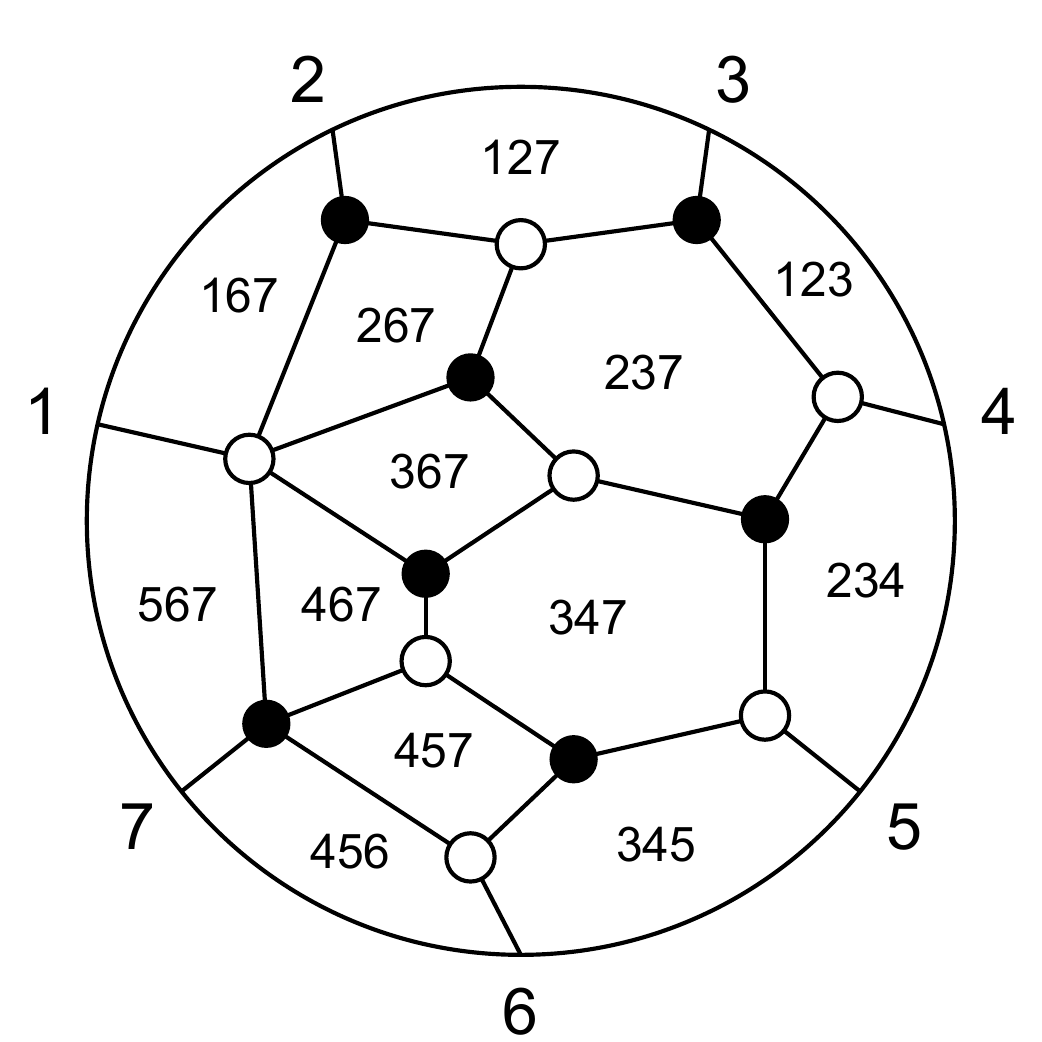}};
		\end{tikzpicture}
		\caption{\label{fig:figure1} A plabic graph with labelled faces.}
	\end{figure}
	
	\begin{definition}\label{definition trip}
		A \textbf{trip} $T$ in a plabic graph $G$ is a path only containing inner edges that turns maximally left at any white vertex, and maximally right at any black vertex.
	\end{definition}
	
	\begin{figure}[!htb]
		\centering
		\begin{tikzpicture}
			\node (r1) at (-2.5, 0) {\includegraphics[height = 2cm]{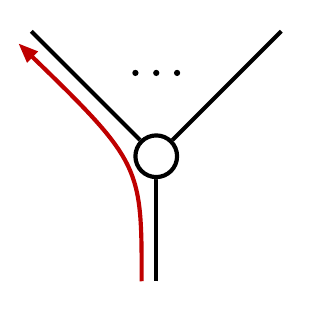}};
			\node (r2) at ( 2.5, 0) {\includegraphics[height = 2cm]{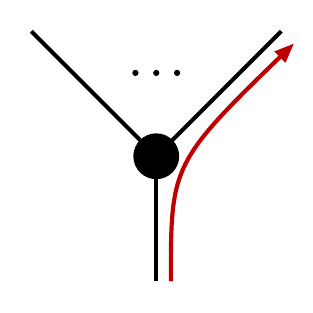}};
		\end{tikzpicture}
		\caption{\label{fig:figure2} Rules of the road: turn left at white vertices and right at black vertices}
	\end{figure}
	
	Trips starting (and thus ending) at the boundary are called \textbf{one way trips}. Any other trip must be a cycle in the interior of $G$ which we refer to as a \textbf{round trip}.
	
	\begin{definition}\label{definition trip-permutation}
		The one way trips inside a plabic graph $G$ determine the \textbf{trip-permutation} $\pi_G \in S_n$ defined by sending $i$ to $j$ if there is a trip from $i$ to $j$ inside $G$.
	\end{definition}
	
	Let $T_1$ and $T_2$ be two (not necessarily distinct) one way trips in some plabic graph. A bicoloured edge passed by both trips in different directions is called an \textbf{essential intersection}. If $T_1$ and $T_2$ have two distinct essential intersections at edges $e_1$ and $e_2$, then we say they have a $\textbf{bad double crossing}$ if both trips traverse first $e_1$ and then $e_2$.
	
	\begin{definition}\label{reducedness}
		A plabic graph $G$ is called \textbf{reduced} if it fulfils the following conditions:
		\begin{itemize}
			\item[\textnormal{(1)}] $G$ contains no round trips.
			\item[\textnormal{(2)}] No trip in $G$ has essential self intersections.
			\item[\textnormal{(3)}] No two trips in $G$ have a bad double crossing.
			\item[\textnormal{(4)}] If $\pi_G(i) = i$, then $i$ is adjacent to a leaf and there are no leaves except at the boundary.
		\end{itemize}
	\end{definition}
	
	\begin{convention}
		All subsequent plabic graphs are assumed to be reduced and admit the trip-permutation $\pi_{k,n} = i \mapsto i+k$.
	\end{convention}

	Under the given convention we know exactly how many inner faces a plabic graph should have (see \cite{Oh2015} Theorem 6.8).
	
	\begin{theorem}\label{number of faces in a plabic graph}
		A plabic graph $G$ contains precisely $k(n-k) + 1$ inner faces.
	\end{theorem}
	
	We will need a way to relate different plabic graphs. On this end we have certain local transformations introduced in the next definition.
	
	\newsavebox{\MovesMIa}
	\sbox{\MovesMIa}{
		\resizebox{2.2cm}{!}{
			\begin{tikzpicture}
				[-,>=stealth',auto, thick]
				\tikzset{white/.style={circle, draw=black, fill=white, scale=1.2}}
				\tikzset{black/.style={circle, draw=black, fill=black, scale=1.2}}
				
				\node[white] (1) at (-0.6, 0.6){};
				\node[black] (2) at ( 0.6, 0.6){};
				\node[white] (3) at ( 0.6,-0.6){};
				\node[black] (4) at (-0.6,-0.6){};
				
				\coordinate (b1) at ($(0,0)+(135:2 and 2)$);
				\coordinate (b2) at ($(0,0)+(45:2 and 2)$);
				\coordinate (b3) at ($(0,0)+(-45:2 and 2)$);
				\coordinate (b4) at ($(0,0)+(-135:2 and 2)$);
				
				\draw
				(1) edge (2)
				(2) edge (3)
				(3) edge (4)
				(4) edge (1)
				
				(1) edge (b1)
				(2) edge (b2)
				(3) edge (b3)
				(4) edge (b4)
				;
	\end{tikzpicture} } }
	
	\newsavebox{\MovesMIb}
	\sbox{\MovesMIb}{
		\resizebox{2.2cm}{!}{
			\begin{tikzpicture}
				[-,>=stealth',auto, thick]
				\tikzset{white/.style={circle, draw=black, fill=white, scale=1.2}}
				\tikzset{black/.style={circle, draw=black, fill=black, scale=1.2}}
				
				\node[black] (1) at (-0.6, 0.6){};
				\node[white] (2) at ( 0.6, 0.6){};
				\node[black] (3) at ( 0.6,-0.6){};
				\node[white] (4) at (-0.6,-0.6){};
				
				\coordinate (b1) at ($(0,0)+(135:2 and 2)$);
				\coordinate (b2) at ($(0,0)+(45:2 and 2)$);
				\coordinate (b3) at ($(0,0)+(-45:2 and 2)$);
				\coordinate (b4) at ($(0,0)+(-135:2 and 2)$);
				
				\draw
				(1) edge (2)
				(2) edge (3)
				(3) edge (4)
				(4) edge (1)
				
				(1) edge (b1)
				(2) edge (b2)
				(3) edge (b3)
				(4) edge (b4)
				;
	\end{tikzpicture} } }
	
	\newsavebox{\MovesMIIa}
	\sbox{\MovesMIIa}{
		\resizebox{2.2cm}{!}{
			\begin{tikzpicture}
				[-,>=stealth',auto, thick]
				\tikzset{white/.style={circle, draw=black, fill=white, scale=1.4}}
				\tikzset{black/.style={circle, draw=black, fill=black, scale=1.4}}
				
				\node[black] (1) at (-0.6, 0){};
				\node[black] (2) at ( 0.6, 0){};
				
				\coordinate (b1) at ($(0,0)+(135:2 and 1.5)$);
				\coordinate (b2) at ($(0,0)+(45:2 and 1.5)$);
				\coordinate (b3) at ($(0,0)+(-45:2 and 1.5)$);
				\coordinate (b4) at ($(0,0)+(-135:2 and 1.5)$);
				\coordinate (b5) at ($(0,0)+(-180:2 and 1.5)$);
				
				\draw
				(1) edge (2)
				
				(1) edge (b1)
				(2) edge (b2)
				(2) edge (b3)
				(1) edge (b4)
				(1) edge (b5)
				;
	\end{tikzpicture} } }
	
	\newsavebox{\MovesMIIb}
	\sbox{\MovesMIIb}{
		\resizebox{2.2cm}{!}{
			\begin{tikzpicture}
				[-,>=stealth',auto, thick]
				\tikzset{white/.style={circle, draw=black, fill=white, scale=1.4}}
				\tikzset{black/.style={circle, draw=black, fill=black, scale=1.4}}
				
				\node[black] (1) at (0, 0){};
				
				\coordinate (b1) at ($(0,0)+(135:2 and 1.5)$);
				\coordinate (b2) at ($(0,0)+(45:2 and 1.5)$);
				\coordinate (b3) at ($(0,0)+(-45:2 and 1.5)$);
				\coordinate (b4) at ($(0,0)+(-135:2 and 1.5)$);
				\coordinate (b5) at ($(0,0)+(-180:2 and 1.5)$);
				
				\draw
				
				(1) edge (b1)
				(1) edge (b2)
				(1) edge (b3)
				(1) edge (b4)
				(1) edge (b5)
				;
	\end{tikzpicture} } }
	
	\newsavebox{\MovesMIIIa}
	\sbox{\MovesMIIIa}{
		\resizebox{2.2cm}{!}{
			\begin{tikzpicture}
				[-,>=stealth',auto, thick]
				\tikzset{white/.style={circle, draw=black, fill=white, scale=0.87}}
				
				\node[white] (1) at (0, 0){};	
				
				\coordinate (b1) at ($(0,0)+(180:1 and 1)$);
				\coordinate (b2) at ($(0,0)+(0:1 and 1)$);
				
				\draw
				(1) edge (b1)
				(1) edge (b2)
				;
	\end{tikzpicture} } }
	
	\newsavebox{\MovesMIIIb}
	\sbox{\MovesMIIIb}{
		\resizebox{2.2cm}{!}{
			\begin{tikzpicture}
				[-,>=stealth',auto, thick]			
				
				\coordinate (b1) at ($(0,0)+(180:1 and 1)$);
				\coordinate (b2) at ($(0,0)+(0:1 and 1)$);
				
				\draw
				(b1) edge (b2)
				;
	\end{tikzpicture} } }
	
	\begin{definition}\label{{definition moves}}
		Let $G$ be plabic graphs then under the correct circumstances we can apply any of the following operations, called \textbf{moves}, to transform $G$ into another plabic graph.
		\begin{itemize}
			\item[(M1)] \textit{Square move:} Suppose $G$ contains a bipartite square face $f$ such that the four inner vertices on $\partial f$ are all of degree $3$. Then we can switch their colours.
			\begin{center}
				\begin{tikzpicture}[align=center]
					\node (A)			{\usebox{\MovesMIa}};
					\node (B) at (5,0)	{\usebox{\MovesMIb} $~~~~~~$};
					
					\draw[to-to, shorten <=2mm, shorten >=2mm]
					(A) edge (B)
					;
				\end{tikzpicture} 
			\end{center}
			
			\item[(M2)] \textit{Unicoloured edge contraction/uncontraction:} We may contract unicoloured edges of $G$ or uncontract inner vertices into unicoloured edges if doing so does not produce any leaves.
			\begin{center}
				\begin{tikzpicture}
					\node (A)			{\usebox{\MovesMIIa}};
					\node (B) at (5,0)	{\usebox{\MovesMIIb} $~~~~~~$};
					
					\draw[to-to, shorten <=2mm, shorten >=2mm]
					(A) edge (B)
					;
				\end{tikzpicture}
			\end{center}
			
			\item[(M3)] \textit{Middle vertex insertion/removal:} Any inner vertex of degree 2 whose neighbours are not both boundary vertices can be removed, merging its edges. On the other hand we may place a vertex of any colour on any inner edge.
			\begin{center}
				\begin{tikzpicture}
					\node (A)			{\usebox{\MovesMIIIa}};
					\node (B) at (5,0)	{\usebox{\MovesMIIIb} $~~~~~~$};
					
					\draw[to-to, shorten <=2mm, shorten >=2mm]
					(A) edge (B)
					;
				\end{tikzpicture}
			\end{center}
		\end{itemize}
	\end{definition}
	
	Note that any bicoloured square face may be transformed via (M1) after applying some of the moves (M2) and (M3). 
	We will usually speak of (M1) as \textbf{mutation} and refer to (M2) and (M3) as \textbf{auxiliary moves}. 
	
	As the name suggests we only really care about plabic graphs up to auxiliary moves. Using the latter we can repeatedly remove degree $2$ vertices and merge unicoloured as long as possible to arrive at a \textbf{contracted} plabic graph. The latter serves as a standard representative. More importantly we have the following theorem (\cite{Postnikov2006} Theorem 13.4).
	
	\begin{theorem}\label{plg_and_wsc_takeaway}
		Any two plabic graphs $G$ and $G'$ can be transformed into another by a sequence of mutations and auxiliary moves. 
	\end{theorem}

	We can label the faces of plabic graphs by subsets of $[n]$. Indeed any trip $T$ in a plabic graph $G$ divides its interior into two parts, one to the left and one to the left right of $T$. 
	
	\begin{definition}\label{face label definition}
		Let $G$ be plabic graph. We define the label of an inner face $f$ of $G$ as the subset of $I(f) \subseteq [n]$ such that $i \in I(f)$ if and only if $f$ lies to the left of the trip starting at $i$.
	\end{definition}

	Using face labels Theorem \ref{plg_and_wsc_takeaway} can be improved: One never needs to mutate a face having the same label in $G$ and $G'$ (see \cite{Oh_2020}). \\
	
	The boundary face incident to $i$ and $i+1$ is always labelled by the set $I_i = \cyc{i-k+1}{i} \subseteq [n]$. We call these labels \textbf{frozen labels}. 
	Note that labels of adjacent inner faces can be obtained from each other by the exchange of one element, since any common edge of them is traversed by two distinct trips in opposite directions. This explains the first part of the following lemma (see \cite{Oh2015} Proposition 8.3.).
	
	\begin{lemma}\label{face label cardinality}
		Let $G$ be plabic graph. Then any inner face of $G$ receives a label of cardinality $k$, and any two labels are distinct.
	\end{lemma}
	
	\newsavebox{\MIvsLabelsA}
	\sbox{\MIvsLabelsA}{
		\resizebox{3.5cm}{!}{
		\begin{tikzpicture}
			[-,>=stealth',auto, thick]
			\tikzset{white/.style={circle, draw=black, fill=white, scale=1.1}}
			\tikzset{black/.style={circle, draw=black, fill=black, scale=1.1}}
			
			\node[white] (1) at (-0.7, 0.7){};
			\node[black] (2) at ( 0.7, 0.7){};
			\node[white] (3) at ( 0.7,-0.7){};
			\node[black] (4) at (-0.7,-0.7){};
			
			\node (Iac) at ( 0  , 0  ) {$Iac$};
			\node (Iad) at (-1.4, 0  ) {$Iad$};
			\node (Iab) at ( 0  , 1.4) {$Iab$};
			\node (Ibc) at ( 1.4, 0  ) {$Ibc$};
			\node (Icd) at ( 0  ,-1.4) {$Icd$};
			
			\coordinate (b1) at ($(0,0)+(135:2.3 and 2.3)$);
			\coordinate (b2) at ($(0,0)+(45:2.3 and 2.3)$);
			\coordinate (b3) at ($(0,0)+(-45:2.3 and 2.3)$);
			\coordinate (b4) at ($(0,0)+(-135:2.3 and 2.3)$);
			
			\draw
			(1) edge (2)
			(2) edge (3)
			(3) edge (4)
			(4) edge (1)
			
			(1) edge (b1)
			(2) edge (b2)
			(3) edge (b3)
			(4) edge (b4)
			;
	\end{tikzpicture} } }
	
	\newsavebox{\MIvsLabelsB}
	\sbox{\MIvsLabelsB}{
		\resizebox{3.5cm}{!}{
		\begin{tikzpicture}
			[-,>=stealth',auto, thick]
			\tikzset{white/.style={circle, draw=black, fill=white, scale=1.1}}
			\tikzset{black/.style={circle, draw=black, fill=black, scale=1.1}}
			
			\node[black] (1) at (-0.7, 0.7){};
			\node[white] (2) at ( 0.7, 0.7){};
			\node[black] (3) at ( 0.7,-0.7){};
			\node[white] (4) at (-0.7,-0.7){};
			
			\node (Ibd) at ( 0  , 0  ) {$Ibd$};
			\node (Iad) at (-1.4, 0  ) {$Iad$};
			\node (Iab) at ( 0  , 1.4) {$Iab$};
			\node (Ibc) at ( 1.4, 0  ) {$Ibc$};
			\node (Icd) at ( 0  ,-1.4) {$Icd$};
			
			\coordinate (b1) at ($(0,0)+(135:2.3 and 2.3)$);
			\coordinate (b2) at ($(0,0)+(45:2.3 and 2.3)$);
			\coordinate (b3) at ($(0,0)+(-45:2.3 and 2.3)$);
			\coordinate (b4) at ($(0,0)+(-135:2.3 and 2.3)$);
			
			\draw
			(1) edge (2)
			(2) edge (3)
			(3) edge (4)
			(4) edge (1)
			
			(1) edge (b1)
			(2) edge (b2)
			(3) edge (b3)
			(4) edge (b4)
			;
	\end{tikzpicture} } }
	
	Mutation is the only move having an effect on labels. This is illustrated below.
	
	\begin{center}
		\begin{tikzpicture}
			\node 			(A) {\usebox{\MIvsLabelsA}};
			\node at (6,0)	(B) {\usebox{\MIvsLabelsB}};
			
			\draw[to-to, shorten <=2mm, shorten >=2mm]
			(A) edge node[above] {(M1)} (B)
			;
			
		\end{tikzpicture}
	\end{center}

	The face labels of a plabic graph come with a nice structure.
	
	\begin{definition}\label{definition weak separation}
		Let $I, J$ be $k$-subsets of $[n]$, then we call $I$ and $J$ \textbf{weakly separated} if we cannot find elements $a, c \in I \setminus J$ and $b, d \in J \setminus I$ such that $(a, b, c, d)$ is strictly cyclically ordered. In this case we write $I \parallel J$. 
		A subset $\mathcal{C} \subseteq \binom{[n]}{k}$ is called a \textbf{weakly separated collection} if its elements are pairwise weakly separated.
	\end{definition}
	
	\begin{theorem}[\cite{Oh2015} Theorem 1.5.]\label{from_plgs_to_wscs}
		If $G$ is a plabic graph then its face labels form a maximal weakly separated collection (with respect to the inclusion of sets) and any such collection comes from a plabic graph in this way. Moreover any plabic graph is uniquely determined by its face labels up to auxiliary moves. 
	\end{theorem}
	
	Using theorem \ref{from_plgs_to_wscs} we will often identify plabic graphs with their face labels and use intuitive language. For example we will call a label in a plabic graph mutable if the corresponding face is and so on.
	
	\begin{remark}\label{looking to the right}
		One can also construct face labels by considering the side right to a trip. Doing this will replace any face label by their complement. Clearly analogues of the results in this section stay true under this alternative convention.
	\end{remark}
	
	\section{Rectangles and Checkboards}
	
	We will now introduce the relevant families of plabic graphs, the first of which is well known and appears for example in \cite{Rietsch_2019} or \cite{FF2018}.
	
	\begin{definition}\label{rectanlge_graph_definition}
		We define the \textbf{rectangle} graph $G^{\text{rec}}_{k,n}$ as follows:
		\begin{itemize}
			\item[(i)]   The inner vertices are arranged on a grid with $2(n-k)$ rows and $2k$ columns, where not all positions are used. Indexing rows and columns by $i$ and $j$, we place an inner vertex $v_{i,j}$ whenever $i+j$ is even, except at positions $(1, 1)$ and $(2(n-k), 2k)$. The vertices $v_{i,j}$ with even $i$ are coloured white, the remaining ones are coloured black.
			
			\item[(ii)]  For any black vertex $v_{i,j}$ add edges to the white vertices $v_{i+1, j \pm 1}$ as long as they exist. Similarly for any white vertex $v_{i,j}$ add edges to the black vertices $v_{i\pm1, j - 1}$ if they exist, and are not already connected to $v_{i,j}$.
			
			\item[(iii)] We place a suitable simple closed curve around the graph constructed so far. The vertices $v_{i,j}$ with $i=1$ or $j=2k$ as well as $v_{2, 2}$ and $v_{2(n-k)-1, 2k-1}$ should be connected to boundary vertices placed on the curve.
			
			\item[(iv)]  Label the boundary vertices on the rectangle clockwise by $1,...,n$, such that the boundary vertex connected to $v_{2, 2}$ receives the label $1$.
		\end{itemize}
		
	\end{definition}

	It is easy to compute the contracted representative of rectangle graphs as defined above. For example figure 
	\ref{fig:figure1} shows a contracted rectangle graph. The representative of definition \ref{rectanlge_graph_definition} is illustrated below.

	\begin{figure}[!htb]
		\centering
		\begin{tikzpicture}
			\node (A) {\includegraphics[height=8cm]{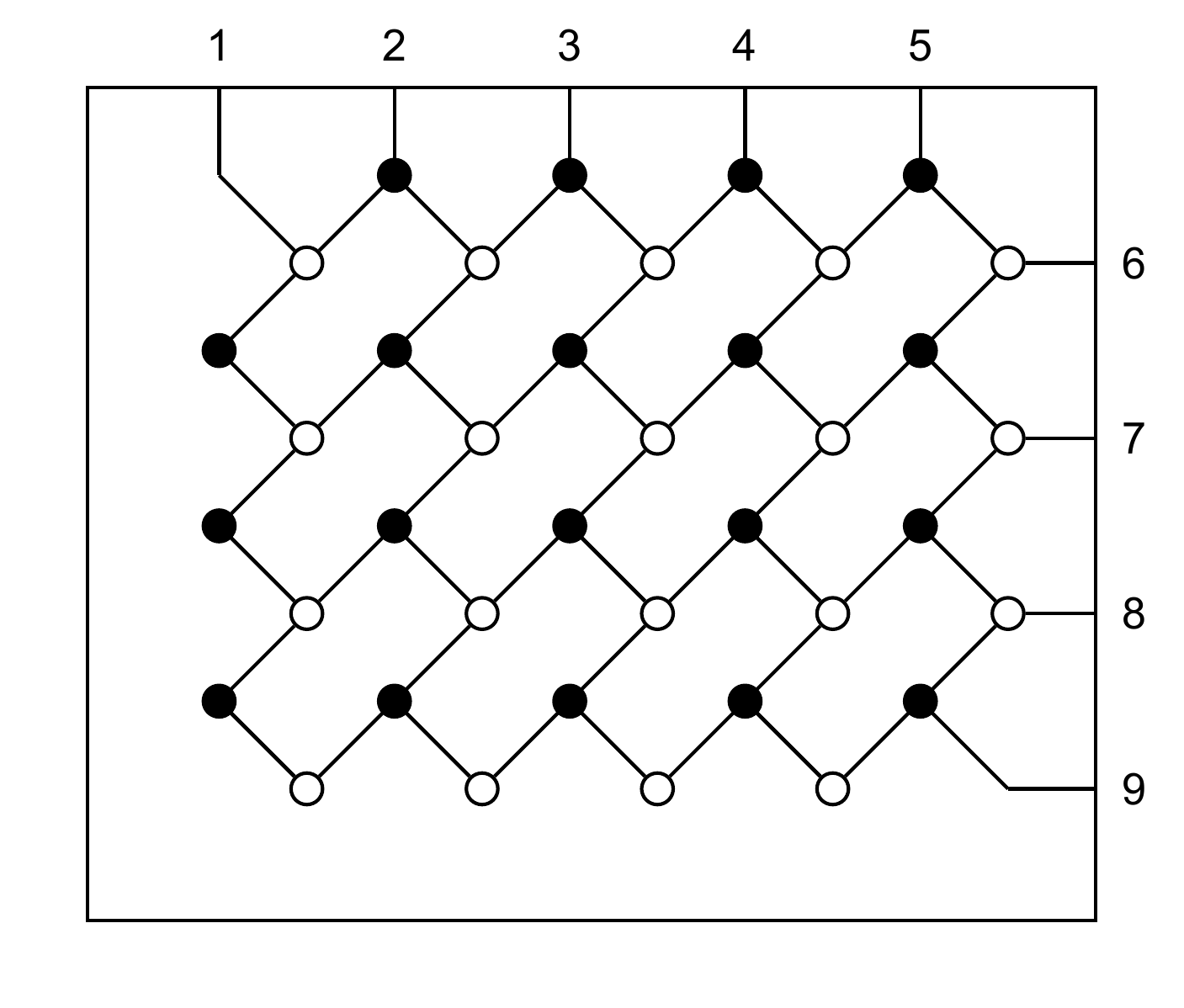} };
			
			\node (0,0) at ($(-3.5, 2.5)$)   {$I_{0,0}$};
			\node (1,0) at ($(-3.5, 1.15)$)  {$I_{1,0}$};
			\node (2,0) at ($(-3.5, -0.2)$)  {$I_{2,0}$};
			\node (3,0) at ($(-3.5, -1.6)$)  {$I_{3,0}$};
			\node (4,0) at ($(-3.5, -2.85)$) {$I_{4,0}$};
			
			\node (0,1) at ($(-2.25, 2.5)$)   {$I_{0,1}$};
			\node (1,1) at ($(-2.25, 1.15)$)  {$I_{1,1}$};
			\node (2,1) at ($(-2.25, -0.2)$)  {$I_{2,1}$};
			\node (3,1) at ($(-2.25, -1.6)$)  {$I_{3,1}$};
			\node (4,1) at ($(-2.25, -2.85)$) {$I_{4,1}$};
			
			\node (0,2) at ($(-0.85, 2.5)$)   {$I_{0,2}$};
			\node (1,2) at ($(-0.85, 1.15)$)  {$I_{1,2}$};
			\node (2,2) at ($(-0.85, -0.2)$)  {$I_{2,2}$};
			\node (3,2) at ($(-0.85, -1.6)$)  {$I_{3,2}$};
			\node (4,2) at ($(-0.85, -2.85)$) {$I_{4,2}$};
			
			\node (0,3) at ($(0.55, 2.5)$)   {$I_{0,3}$};
			\node (1,3) at ($(0.55, 1.15)$)  {$I_{1,3}$};
			\node (2,3) at ($(0.55, -0.2)$)  {$I_{2,3}$};
			\node (3,3) at ($(0.55, -1.6)$)  {$I_{3,3}$};
			\node (4,3) at ($(0.55, -2.85)$) {$I_{4,3}$};
			
			\node (0,4) at ($(1.9, 2.5)$)   {$I_{0,4}$};
			\node (1,4) at ($(1.9, 1.15)$)  {$I_{1,4}$};
			\node (2,4) at ($(1.9, -0.2)$)  {$I_{2,4}$};
			\node (3,4) at ($(1.9, -1.6)$)  {$I_{3,4}$};
			\node (4,4) at ($(1.9, -2.85)$) {$I_{4,4}$};
			
			\node (0,5) at ($(3.25, 2.5)$)   {$I_{0,5}$};
			\node (1,5) at ($(3.25, 1.15)$)  {$I_{1,5}$};
			\node (2,5) at ($(3.25, -0.2)$)  {$I_{2,5}$};
			\node (3,5) at ($(3.25, -1.6)$)  {$I_{3,5}$};
			\node (4,5) at ($(3.25, -2.85)$) {$I_{4,5}$};
		\end{tikzpicture}
		\caption{\label{fig:figure3} The rectangle graph $G^\text{rec}_{5,9}$ with its face labels $I_{i,j}$ arranged in a grid.}
	\end{figure}
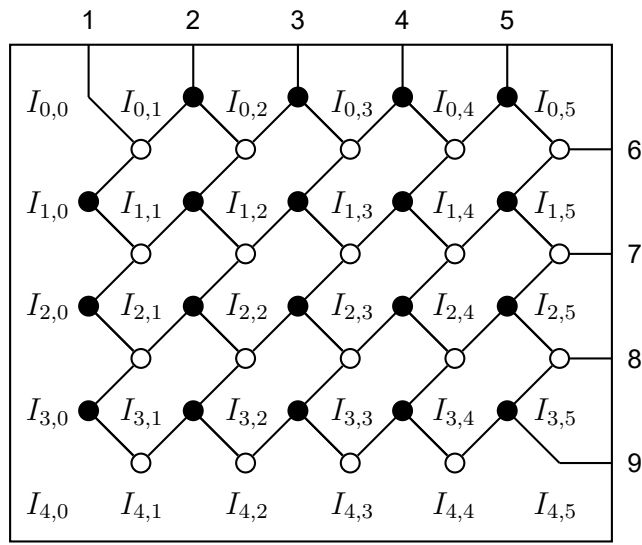
	
	The face labels of $G^{\text{rec}}_{k,n}$ are well known. 
	
	\begin{proposition}\label{rectangle_graph_labels}
		Let us arrange the faces of $G^{\text{rec}}_{k,n}$ being not adjacent to its boundary on a grid with $n-k-1$ rows and $k-1$ columns. Indexing rows by $i$ and columns by $j$, they are labelled by 
		$ I^\text{rec}_{i,j} = L^\text{rec}_{i,j} \cup R^\text{rec}_{i,j} $, where
		$$L^\text{rec}_{i,j} = [i+j] \setminus [i] ~ \text{ and } ~ R^\text{rec}_{i,j} = [n] \setminus [n-k+j].$$
	\end{proposition}
	
	\begin{remark}\label{extended rectangle grid}
		Allowing $i = 0, n-k$ and $j = 0,k$ in the formulas of \ref{rectangle_graph_labels} results in the \enquote{correct} frozen labels as indicated in figure \ref{fig:figure2}.
	\end{remark}
	
	Let us now introduce the second family of graphs, which is the one we are primarily concerned with. This family of graphs is inspired by the Quadrilateral Arrangements in \cite{Scott2006}. These objects were used by Scott to show that the coordinate rings of Grassmannians carry the structure of a cluster algebra.
	
	\begin{definition}\label{checkboard graphs}
		We construct the \textbf{checkboard} graph $G^{\text{ch}}_{k,n}$ by the following steps:
		\begin{itemize}
			\item[(i)]   The inner vertices are arranged on a grid with $n-k$ rows and $k$ columns, where all indices are used. Let them be denoted by $v_{i,j}$. The colour of $v_{i,j}$ should be white whenever $i+j$ is even and black otherwise.
			
			\item[(ii)]  Connect any inner vertices $v_{i,j}$ and $v_{l,m}$ whose indices satisfy $|i-l|+|j-m|=1$.
			
			\item[(iii)] We place a suitable simple closed curve around the graph constructed so far. The white vertices at the left and right side of the grid, as well as the black vertices at its top and bottom should be connected to boundary vertices lying on the curve.
			
			\item[(iv)]  Label the boundary vertices on the rectangle clockwise by $1,...,n$, such that the boundary vertex connected to the white top left vertex $v_{1,1}$ receives the label $n$.
		\end{itemize}
	\end{definition}
	
	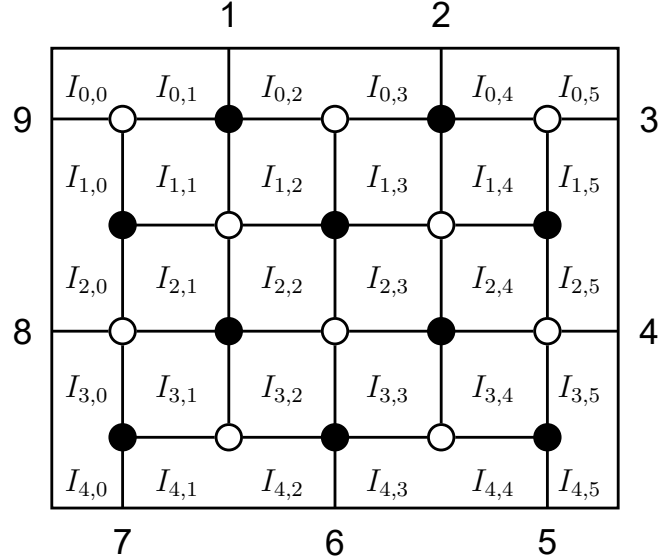
\begin{figure}[!htb]
		\centering
		\begin{tikzpicture}
			\node (A) {\includegraphics[height=8cm]{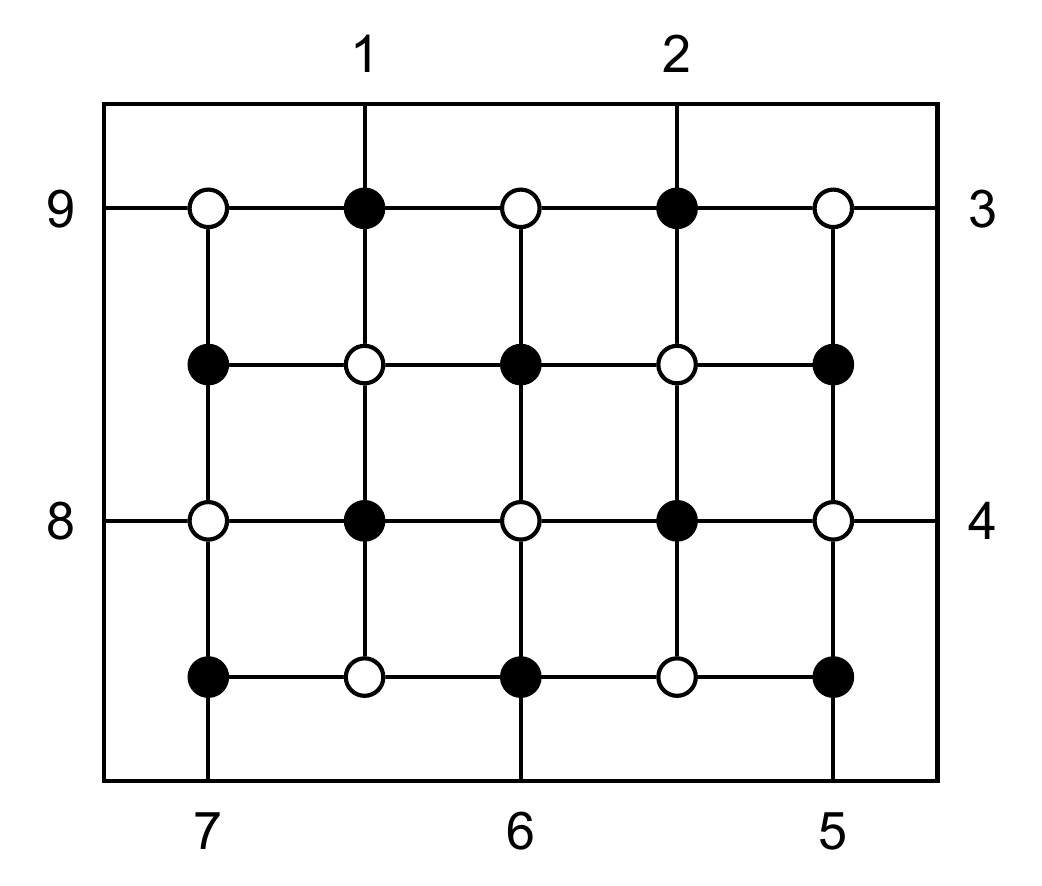} };
			
			\node (0,0) at ($(-3.3, 2.5)$)   {$I_{0,0}$};
			\node (1,0) at ($(-3.3, 1.3)$)  {$I_{1,0}$};
			\node (2,0) at ($(-3.3, -0.05)$)  {$I_{2,0}$};
			\node (3,0) at ($(-3.3, -1.5)$)  {$I_{3,0}$};
			\node (4,0) at ($(-3.3, -2.75)$) {$I_{4,0}$};
			
			\node (0,1) at ($(-2.1, 2.5)$)   {$I_{0,1}$};
			\node (1,1) at ($(-2.1, 1.3)$)  {$I_{1,1}$};
			\node (2,1) at ($(-2.1, -0.05)$)  {$I_{2,1}$};
			\node (3,1) at ($(-2.1, -1.5)$)  {$I_{3,1}$};
			\node (4,1) at ($(-2.1, -2.75)$) {$I_{4,1}$};
			
			\node (0,2) at ($(-0.7, 2.5)$)   {$I_{0,2}$};
			\node (1,2) at ($(-0.7, 1.3)$)  {$I_{1,2}$};
			\node (2,2) at ($(-0.7, -0.05)$)  {$I_{2,2}$};
			\node (3,2) at ($(-0.7, -1.5)$)  {$I_{3,2}$};
			\node (4,2) at ($(-0.7, -2.75)$) {$I_{4,2}$};
			
			\node (0,3) at ($(0.7, 2.5)$)   {$I_{0,3}$};
			\node (1,3) at ($(0.7, 1.3)$)  {$I_{1,3}$};
			\node (2,3) at ($(0.7, -0.05)$)  {$I_{2,3}$};
			\node (3,3) at ($(0.7, -1.5)$)  {$I_{3,3}$};
			\node (4,3) at ($(0.7, -2.75)$) {$I_{4,3}$};
			
			\node (0,4) at ($(2.1, 2.5)$)   {$I_{0,4}$};
			\node (1,4) at ($(2.1, 1.3)$)  {$I_{1,4}$};
			\node (2,4) at ($(2.1, -0.05)$)  {$I_{2,4}$};
			\node (3,4) at ($(2.1, -1.5)$)  {$I_{3,4}$};
			\node (4,4) at ($(2.1, -2.75)$) {$I_{4,4}$};
			
			\node (0,5) at ($(3.25, 2.5)$)   {$I_{0,5}$};
			\node (1,5) at ($(3.25, 1.3)$)  {$I_{1,5}$};
			\node (2,5) at ($(3.25, -0.05)$)  {$I_{2,5}$};
			\node (3,5) at ($(3.25, -1.5)$)  {$I_{3,5}$};
			\node (4,5) at ($(3.25, -2.75)$) {$I_{4,5}$};
		\end{tikzpicture}
		\caption{\label{fig:figure4} The checkboard graph $G^\text{ch}_{5,9}$ with its face labels $I_{i,j}$ arranged in a grid.}
	\end{figure}
	
	We also know the labels of checkboard graphs.
	\begin{proposition}\label{checkboard_graph_labels}
		Let us arrange the faces of $G^{\text{ch}}_{k,n}$ not adjacent to the boundary on a grid with $n-k-1$ rows and $k-1$ columns. Indexing rows by $i$ and columns by $j$, they are labelled by 
		$ I^\text{ch}_{i,j} = L^\text{ch}_{i,j} \cup R^\text{ch}_{i,j} $, where
		$$L^\text{ch}_{i,j} = \sigma^{- \lceil (i+j)/2 \rceil}([i+j] \setminus [i]) ~ \text{ and } ~ R^\text{ch}_{i,j} = \sigma^{- \lceil (i+j)/2 \rceil}([n] \setminus [n-k+j]). $$
		Here $\sigma = i \mapsto i+1$ denotes the clockwise rotation by one step.
	\end{proposition}
 	\begin{proof}
 		The proof boils down to computing which trips pass the vertical edges of $G^{\text{ch}}_{k,n}$ in each row, which can be done inductively.
 		Then the label of a face depends only on the trips passing the vertical edges in the same row. We leave the details to the reader.
 	\end{proof}
	
	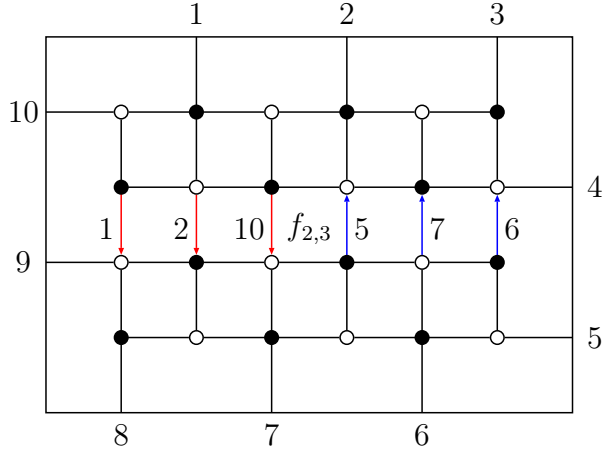
\begin{figure}[!htb]
		\centering
		\begin{tikzpicture}
			\node (background) at ($(0, 0)$){\includegraphics[width = 8cm, height = 6 cm]{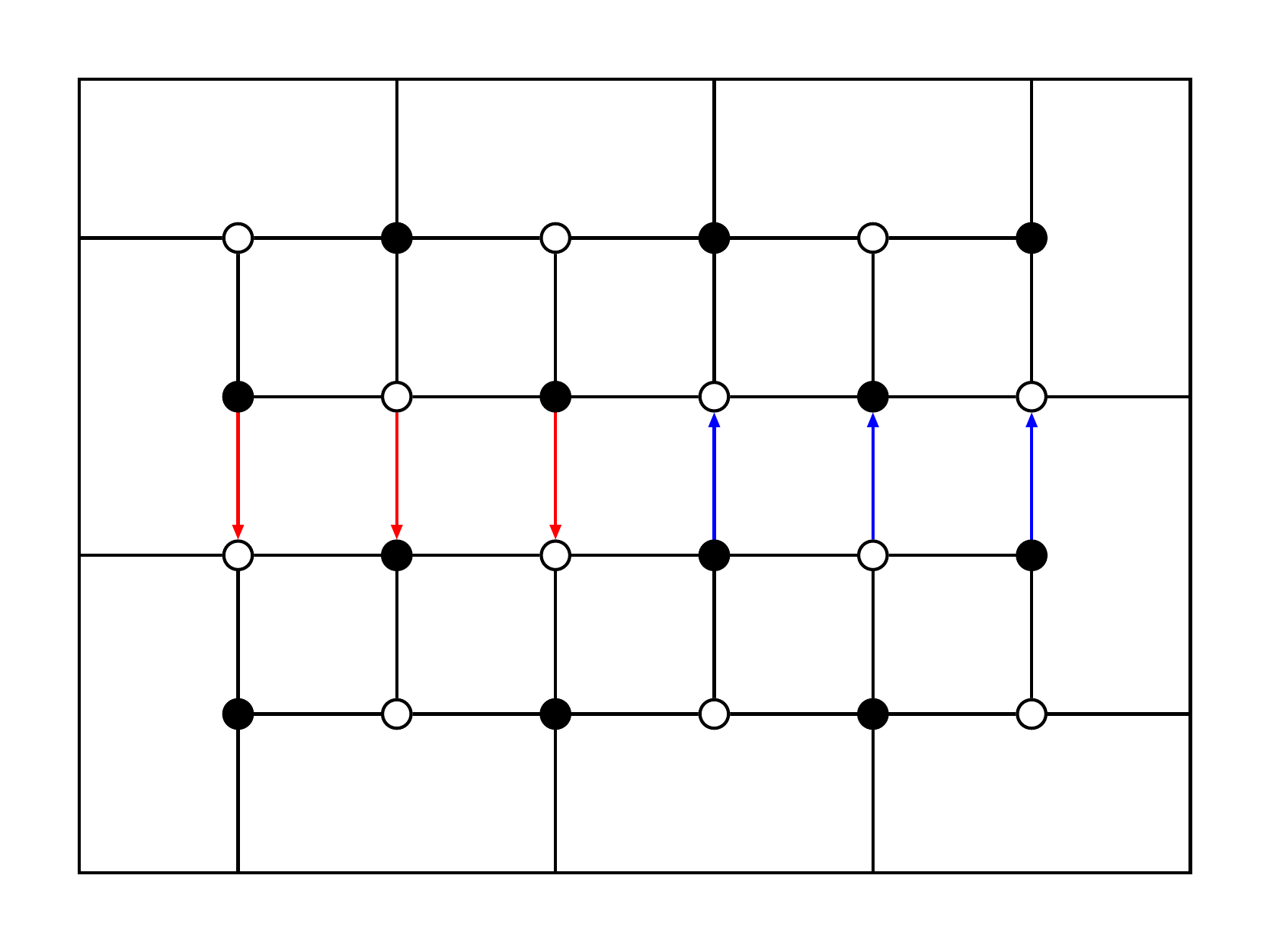}};
			
			\node (b1)  at ($( -1.5,  2.8 )$) {$1$};
			\node (b2)  at ($(  0.5,  2.8 )$) {$2$};
			\node (b3)  at ($(  2.5,  2.8 )$) {$3$};
			\node (b4)  at ($(  3.8,  0.5 )$) {$4$};
			\node (b5)  at ($(  3.8, -1.5 )$) {$5$};
			\node (b6)  at ($(  1.5, -2.8 )$) {$6$};
			\node (b7)  at ($( -0.5, -2.8 )$) {$7$};
			\node (b8)  at ($( -2.5, -2.8 )$) {$8$};
			\node (b9)  at ($( -3.8, -0.5 )$) {$9$};
			\node (b10) at ($( -3.8,  1.5 )$) {$10$};
			
			\node (fij) at ($(0, -0.05)$) {$f_{2,3}$};
			
			\node (d1) at ($(-0.8, -0.05)$) {$10$};
			\node (d2) at ($(-1.7, -0.05)$) {$2$};
			\node (d3) at ($(-2.7, -0.05)$) {$1$};
			
			\node (u1) at ($( 0.7, -0.05)$) {$5$};
			\node (u2) at ($( 1.7, -0.05)$) {$7$};
			\node (u3) at ($( 2.7, -0.05)$) {$6$};
		\end{tikzpicture}
		\caption{\label{fig:figure5} The edges contributing to the label of the face $f_{2,3}$ in row $2$ and column $3$ of the checkboard graph $G^\text{ch}_{6,10}$.}
	\end{figure}
	It is easy to see along the way of computing labels that the trip permutation of $G^\text{ch}_{k,n}$ is in fact $\pi_{k,n}$. To obtain reducedness one can either check Definition \ref{reducedness} directly or use the criterion found in chapter 10.2 of \cite{KW2011}.
	
	\begin{remark}\label{extended checkboard grid}
		Allowing $i = 0, n-k$ and $j = 0,k$ in the formulas of \ref{checkboard_graph_labels} gives the \enquote{correct} frozen labels as indicated in figure \ref{fig:figure3}
	\end{remark}
	
	Before describing two more families of graphs we need to make a quick break and describe what happens when taking complements of labels (see also \cite{FF2018}).
	
	\begin{lemma}\label{complement_of_wsc}
		Let $\mathcal{C} \subseteq \binom{[n]}{k}$ be a maximal weakly separated collection. Then $\mathcal{C}' = \{ [n] \setminus I \mid I \in \mathcal{C} \}$ is a maximal weakly separated collection in $\binom{[n]}{n-k}$.
	\end{lemma}
	
	\begin{proof}
		Let $G$ be the plabic graph corresponding to $\mathcal{C}$. Consider what effect inverting the colours of vertices in $G$ has on its trips. The trip from $i$ to $i+k$ will be reversed to a trip from $i+k$ to $i$. Thus we obtain a plabic graph $G''$ with trip-permutation $\pi_{n-k,n}$. Now define $G'$ by replacing the label of each boundary vertex $i$ in $G''$ with $i-k$. Observe that the face labels of $G'$ are precisely given by $\mathcal{C}'$.
	\end{proof}
	
	\begin{definition}
		The \textbf{dual rectangle graph} $\widehat G^{\text{rec}}_{k,n}$ is defined to be the plabic graph whose face labels are obtained by taking complements of the face labels of $G^{\text{rec}}_{n-k,n}$. 
		Similarly we define the \textbf{dual checkboard graph} $\widehat G^{\text{ch}}_{k,n}$.
	\end{definition}
	
	Note that in contrast to the definition in \cite{FF2018}, $\widehat G^{\text{rec}}_{k,n}$ and $\widehat G^{\text{ch}}_{k,n}$ are again plabic graphs with trip-permutation $\pi_{k,n}$ as they are computed from plabic graphs with trip-permutation $\pi_{n-k,n}$. Thus they are related by sequences of mutations.
	The proof of lemma \ref{complement_of_wsc} also gives us a recipe for constructing these graphs explicitly. For $\widehat G^{\text{rec}}_{k,n}$ just invert the colours of $G^{\text{rec}}_{n-k,n}$ and replace boundary vertex $i$ with $i- (n-k) \equiv i+k$. 
	
	\begin{proposition}
		The face labels of $\widehat G^{\text{rec}}_{k,n}$ arranged on a grid with $k-1$ rows and $n-k-1$ columns are given by $\widehat I^\text{rec}_{i,j} = \widehat L^\text{rec}_{i,j} \cup \widehat R^\text{rec}_{i,j}$, where
		$$ \widehat L^\text{rec}_{i,j} = [i] ~ \text{ and } ~ \widehat R^\text{rec}_{i,j} = [k+j] \setminus [i+j] .$$
		Similarly the face labels of $\widehat G^{\text{ch}}_{k,n}$ are given by $\widehat I^\text{ch}_{i,j} = \widehat L^\text{ch}_{i,j} \cup \widehat R^\text{ch}_{i,j}$, where
		$$ \widehat L^\text{ch}_{i,j} = \sigma^{- \lceil (i+j)/2 \rceil}([i]) ~ \text{ and } ~ \widehat R^\text{ch}_{i,j} = \sigma^{- \lceil (i+j)/2 \rceil}([k+j] \setminus [i+j]) .$$
		Here again $\sigma = i \mapsto i+1$.
	\end{proposition}

	\section{Dihedral symmetry}
	We explicitly describe a natural action of the dihedral group.
	Let $\sigma = i \mapsto i+1$ and $\tau = i \mapsto n+2-i$ 
	be permutations in $S_n$. We realize the dihedral group $D_n$ with $2n$ elements, as the subgroup of $S_n$ generated by $\sigma$ and $\tau$. Recall Definition \ref{definition weak separation} of weak separation.
	
	We say a permutation $\rho \in S_n$ \textbf{preserves weak separation} if $I \parallel J$ implies $\rho(I) \parallel \rho(J)$ for any two $k$-subsets $I, J \subseteq [n]$. This property obviously defines a subgroup of $S_n$. 
	
	\begin{proposition}\label{D_n_preserves_ws}
		Let $\rho \in S_n$. Then $\rho$ preserves weak separation if and only if $\rho \in D_n$. In particular $D_n$ acts on the set of maximal weakly separated collections.
	\end{proposition}
	\begin{proof}
		It is obvious that any $\rho \in D_n$ preserves weak separation, since it either preserves or reverses cyclic orders. Now suppose $\rho \in S_n$ preserves weak separation. Possibly after replacing $\rho$ with $\sigma^m \rho$ for some $m \in \Z$, we may assume $\rho(1) = 1$. \\
		Now choose any $1 < a < b < c \leq n$. Since we require $2 \leq k \leq n-2$ there are at least $k-2$ elements remaining in $[n] \setminus \{1,a,b,c\}$. Thus we can extend $\{1,a\}$ and $\{b,c\}$ by a common $k-2$-subset of $[n]$ to obtain $I$ and $J$ such that $I \setminus J = \{1,a\}$ and $J \setminus I = \{b,c\}$. By construction we then have $I \parallel J$, hence $\rho(I) \parallel \rho(J)$. \\
		Notice that we cannot have $\rho(b) < \rho(a) < \rho(c)$ or $\rho(c) < \rho(a) < \rho(b)$ as this would contradict $\rho(I) \parallel \rho(J)$. Hence either $\rho(a) < \rho(b), \rho(c)$ or $\rho(b), \rho(c) < \rho(a)$. \\
		Assume the former and take any $d > a$. If we had $\rho(d) < \rho(a)$, then either $1 < a < d < c$ and $\rho(d) < \rho(a) < \rho(c)$ or $1 < a  < b < d$ and $\rho(d) < \rho(a) < \rho(b)$. Both possibilities contradict the previously established observation, so we conclude that $\rho(a) < \rho(d)$. Similarly the other case leads to $\rho(d) < \rho(a)$ for any $d > a$. \\
		This proves: If $2 \leq a \leq n-2$ then either $\rho(a) < \rho(d)$ for all $d > a$ or $\rho(d) < \rho(a)$ for all $d > a$. \\
		Now consider $a = 2$. If $\rho(d) > \rho(a)$ for all $d > a$ then we must have $\rho(2) = 2$. Otherwise we obtain $\rho(2) = n$. We can assume $\rho(2) = 2$ by replacing $\rho$ with $\tau\rho$ if needed. \\
		Suppose $\rho \neq \text{id}$ and let $i \in [n]$ be its minimal non fixpoint. Note that $2 < i < n$, If $i = n-1$ then we immediately obtain $\rho(i) = n$. Otherwise we must have $\rho(i) > \rho(d)$ for all $d > i$, so that $\rho(i) = n$ follows in any case.
		To derive a contradiction we construct $k$-subsets $I,J$ such that $I \setminus J = \{1,i\}$ and $J \setminus I = \{2,i+1\}$. Clearly $I$ and $J$ are not weakly separated, but $\rho(I)$ and $\rho(J)$ are. This is a contradiction as $\rho^{-1}$ must also preserve weak separation. Thus $\rho = \text{id}$. 
	\end{proof}
	
	\begin{figure}[!htb]
		\centering
		\begin{tikzpicture}
			\node (1) at (0,0) {\includegraphics[height=5cm]{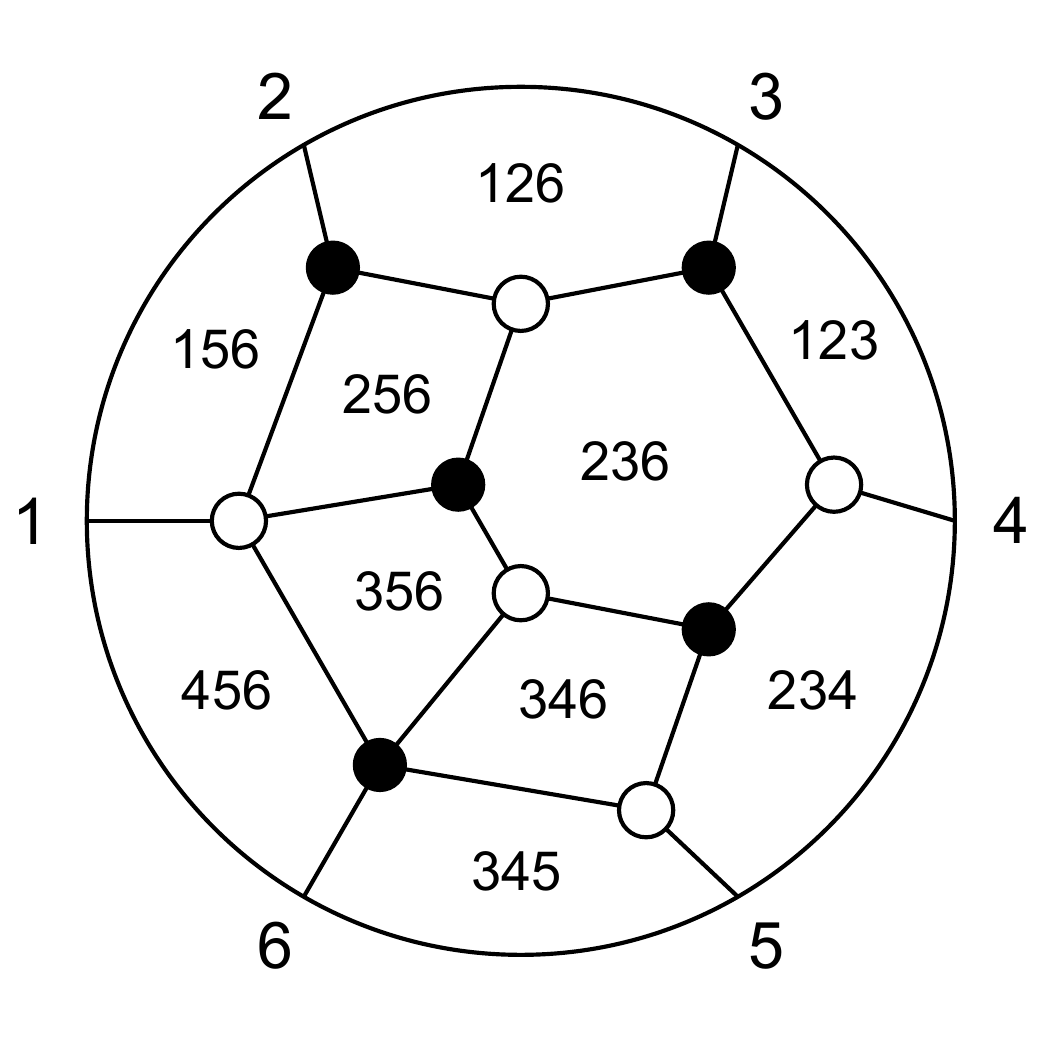}};
			\node (2) at (8.5,0) {\includegraphics[height=5cm]{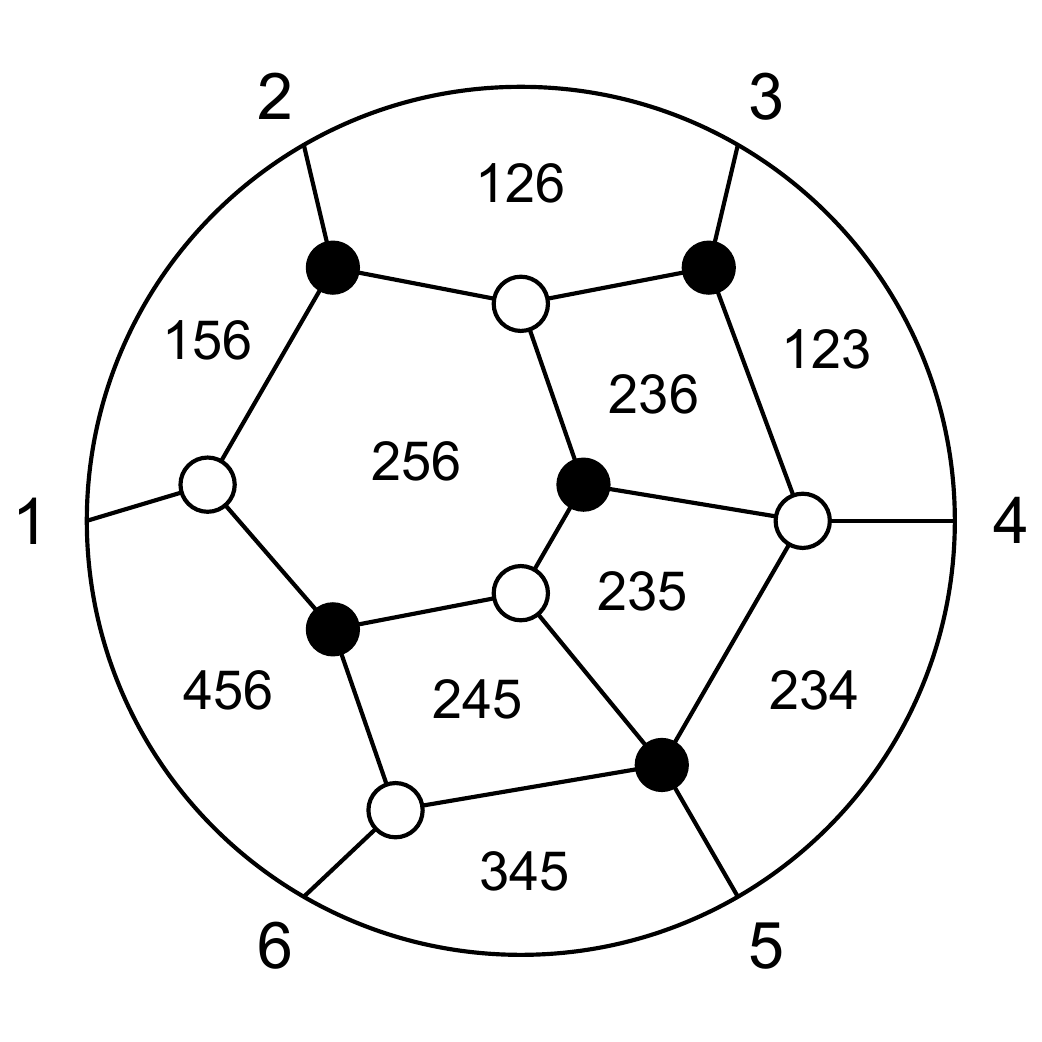}};
			
			\draw[-to, shorten >= 2mm, shorten <= 2mm]
			(1) edge node[above] {$(2~6)(3~5)$} (2)
			;
		\end{tikzpicture}
		\caption{\label{fig:figure6} Reflecting a plabic graph ($k = 3, n = 6$).}
	\end{figure}
	
	We obtain an action on plabic graphs by identifying them with their labels. The following lemma describes this explicitly.
	
	\begin{lemma}\label{action on plgs}
		Let $G$ be a plabic graph of type $\pi_{k,n}$ and $\rho \in D_n$. If $\rho$ is a rotation then $\rho G$ can be obtained from $G$ by replacing the boundary vertex $i$ with $\rho(i)$. Otherwise if $\rho$ is a reflection, then we obtain $\rho G$ as reflection of $G$ at some line together with the relabelling $i \mapsto \sigma^k\rho(i)$. 
	\end{lemma}
	
	For later use we observe that mutation is compatible with this action in the obvious way.
	
	\begin{lemma}\label{mutation_vs_action}
		Let $\rho \in D_n$. If $G$ a plabic graph in which $I$ labels a square face, then $\rho(I)$ labels a square face in $\rho G$ and we have
		$$ \rho \mu_I G = \mu_{\rho(I)} \rho G,$$ where $\mu_I$ denotes the mutation of the face labeled by $I$.
	\end{lemma}

	Let us examine orbits and stabilizers. If $k = 2, n = 4$ then there are only two plabic graphs (up to auxiliary moves). Apart from boundary faces they both contain precisely one square face and arise as a rotation from another. Let us now assume $n > 4$.
	
	\begin{proposition}\label{dual_rec_is_reflection_of_rec}
		The stabilizer of $G^{\text{rec}}_{k,n}$ is trivial whenever $n \geq 5$.
		Moreover we always have $$\widehat G^{\text{rec}}_{k,n} = \rho G^{\text{rec}}_{k,n}, $$ where $\rho = \tau \sigma = i \mapsto n+1-i$.
	\end{proposition}

	\begin{proof}
		We compare $\widehat G^{\text{rec}}_{k,n}$ and $\rho G^{\text{rec}}_{k,n}$ by their labels.
		For $i = 1, \ldots, n-k-1$ and $j = 1, \ldots, k-1$ we compute
		$$ \rho(L^\text{rec}_{i,j}) = \widehat R^\text{rec}_{k-j, n-k-i} 
		~~\text{ and }~~ 
		\rho(R^\text{rec}_{i,j}) = \widehat L^\text{rec}_{k-j, n-k-i}. $$
		Thus we have $\rho(I^\text{rec}_{i,j}) = \widehat I^\text{rec}_{k-j,n-k-i}$. The frozen labels are obviously send onto frozen labels, so we conclude $\widehat G^{\text{rec}}_{k,n} = \rho G^{\text{rec}}$. \\
		As for the stabilizers suppose $n \geq 5$. Then in the contracted representative of $G^\text{rec}_{k,n}$, the boundary vertex $b_1$ is adjacent to a white vertex, $b_2, \ldots, b_k$ are adjacent to black vertices, $b_{k+1}, \ldots, b_{n-1}$ are adjacent to white vertices and $b_n$ is adjacent to a black vertex again. Since either $k-1 > 1$ or $n-k-1 > 1$, it is clear that no rotation or reflection can preserve this property. 
	\end{proof}
	
	Sadly, life is not that easy in case of the checkboard graphs.
	In general there are elements of $D_n$ that fix $G^{\text{ch}}_{k,n}$ even for non trivial $k, n$ and $\widehat G^{\text{ch}}_{k,n}$ in general does not lie in the orbit of $G^{\text{ch}}_{k,n}$.
	
	\begin{proposition}\label{check vs dual conjugate}
		Whenever both $k$ and $n-k$ are odd, the orbit of $G^{\text{ch}}_{k,n}$ does not contain $\widehat G^{\text{ch}}_{k,n}$. Otherwise we have 
		$$ \widehat G^{\text{ch}}_{k,n} = \sigma^m G^{\text{ch}}_{k,n} ,$$
		where $m = k/2$ if $k$ is even and $m = -(n-k)/2$ if $n-k$ is even.
	\end{proposition}
	
	\begin{proof}
		If $k$ and $n-k$ are odd, then all four corners in the grid of $G^\text{ch}_{k,n}$ are white, while they are black in $\widehat G^\text{ch}_{k,n}$. Thus the latter cannot be a rotation or reflection of $G^\text{ch}_{k,n}$. \\
		Now for $i = 1, \ldots, n-k-1$ and $j = 1, \ldots, k-1$ one can check
		$$ \sigma^{k/2} L^\text{ch}_{i,j} = \widehat R^\text{ch}_{k-j,i} ~ \text{ and } ~ \sigma^{k/2} R^\text{ch}_{i,j} = \widehat L^\text{ch}_{k-j,i} ,$$ if $k$ is even. 
		Similar if $n-k$ is even the identities
		$$ \sigma^{-(n-k)/2} L^\text{ch}_{i,j} = \widehat R^\text{ch}_{j,n-k-i} ~ \text{ and } ~ \sigma^{-(n-k)/2} R^\text{ch}_{i,j} = \widehat L^\text{ch}_{j,n-k-i}$$ 
		hold.
		In any case, labels of $G^\text{ch}_{k,n}$ are send to labels of $\widehat G^{\text{ch}}_{k,n}$ and we are done.
	\end{proof}
	
	Of course we can also calculate the stabilizers of the (dual) checkboards.
	
	\begin{proposition}\label{stabilizers of checkboards}
		Let $S$ be the stabilizer of $G^{\text{ch}}_{k,n}$ and $\hat S $ be the one of $\widehat G^{\text{ch}}_{k,n}$. Then we have
		$$ S = 
		\begin{cases}
			\langle \sigma^{n/2} \rangle 			&2| k, 2 | n \\
			\langle \rho, \sigma^{n/2} \rangle		&2 \nmid k, 2 | n \\
			\langle \rho \rangle 					&2 \nmid n
		\end{cases}
		~~~\text{and}~~~
		\hat S = 
		\begin{cases}
			\langle \sigma^{n/2} \rangle 					&2 | k, 2 | n \\
			\langle \sigma^k \rho, \sigma^{n/2} \rangle 	&2 \nmid k, 2 | n \\
			\langle \sigma^k \rho \rangle 					&2 \nmid n 
		\end{cases}
		$$
		
		where $\rho = x \mapsto \frac{(k-1)(n-1)}{2} - x$, is defined whenever $(k-1)(n-1)$ is even.
	\end{proposition}
	
	\begin{proof} This easily follows from lemma \ref{action on plgs} and by considering the labels of the boundary vertices adjacent to the four corners of the \enquote{grids} inside $G^{\text{ch}}_{k,n}$ resp. $\widehat G^{\text{ch}}_{k,n}$.
	\end{proof}
	
	An easy consequence of the previous two propositions is that we always have $$|D_n G^\text{ch}_{k,n} \cup D_n \widehat G^\text{ch}_{k,n}| = n.$$
	
	\section{Quivers}
	
	We briefly describe quivers and their mutations, as well as how to associate quivers to plabic graphs. Here we follow mostly \cite{Rietsch_2019}, but also want to mention \cite{Fomin2016}.
	
	\begin{definition}\label{definition quiver}
		A \textbf{quiver} $Q$ is a directed graph. We shall always assume that $Q$ has no loops or two cycles, but allow multi edges. Moreover each vertex is designated as being either \textbf{frozen} or \textbf{mutable}.
	\end{definition}
	
	\begin{definition}\label{quiver mutation}
		Let $Q$ be quiver and $w \in Q$ a mutable vertex. We define the \textbf{quiver mutation} in direction $w$ which transforms $Q$ into another quiver $Q' = \mu_w(Q)$ by the following steps:
		
		\begin{itemize}
			\item[(1)] For every oriented path $i \to w \to j$, add a new arrow $i \to j$ if at least one of $i$ and $j$ is mutable.
			\item[(2)] Reverse all arrows incident to $w$.
			\item[(3)] Repeatedly remove $2$-cycles until none are left.
		\end{itemize}
		
	\end{definition}
	
	\newsavebox{\quiverMutationExI}
	\sbox{\quiverMutationExI}{
		\begin{tikzpicture}
			\tikzstyle{every path}=[line width=1pt, shorten <=1mm, shorten >=1mm]
			\tikzset{black/.style={circle, draw=black, fill=black, scale=0.5}}
			\tikzset{blue/.style={circle, draw=blue, fill=blue, scale=0.5}}
			
			\node[blue]  (A) at (-0.2 , 0  ) {};
			\node[black] (B) at ( 1   ,-0.2) {};
			\node[black] (C) at ( 2.1 ,-0.5) {};
			\node (label) at ( 2.1 ,-0.2) {$w$};
			\node[blue]  (D) at (2.6  ,-1.5) {};
			\node[blue]  (E) at ( 2   , 0.5) {};
			\node[black] (F) at ( 3.2 ,-0.3) {};
			
			\draw
			[-latex](A) edge (B)
			[-latex](B) edge[bend left=15] (C)
			[-latex](B) edge[bend right=15] (C)
			[-latex](C) edge (D)
			[-latex](D) edge (F)
			[-latex](F) edge (C)
			[-latex](B) edge (E)
			;
	\end{tikzpicture} }
	
	\newsavebox{\quiverMutationExII}
	\sbox{\quiverMutationExII}{
		\begin{tikzpicture}
			\tikzstyle{every path}=[line width=1pt, shorten <=1mm, shorten >=1mm]
			\tikzset{black/.style={circle, draw=black, fill=black, scale=0.5}}
			\tikzset{blue/.style={circle, draw=blue, fill=blue, scale=0.5}}
			
			\node[blue]  (A) at (-0.2 , 0  ) {};
			\node[black] (B) at ( 1   ,-0.2) {};
			\node[black] (C) at ( 2.1 ,-0.5) {};
			\node (label) at ( 2.1 ,-0.2) {$w$};
			\node[blue]  (D) at (2.6  ,-1.5) {};
			\node[blue]  (E) at ( 2   , 0.5) {};
			\node[black] (F) at ( 3.2 ,-0.3) {};
			
			\draw
			[-latex](A) edge (B)
			[-latex](C) edge[bend left=15]  (B)
			[-latex](C) edge[bend right=15] (B)
			[-latex](D) edge (C)
			[-latex](C) edge (F)
			[-latex](B) edge (E)
			[-latex](B) edge[bend right = 20] (D)
			[-latex](B) edge[bend right = 50] (D)
			;
	\end{tikzpicture} }
	
	\begin{figure}[!htb]
		\centering
		\begin{tikzpicture}
			\tikzstyle{every node}=[text width={width("$~~~~~~~~~~~~~~~~~~~~~~~~~~~~~~~~~~~~~~~~~~~$")},align=center]
			\node			  (A) {\usebox{\quiverMutationExI}};
			\node[right=of A] (B) {\usebox{\quiverMutationExII}};
			
			\draw
			[-to](A) edge node[above] {$\mu_w$} (B)
			;
		\end{tikzpicture}
		\caption{\label{fig:figure7} A quiver mutation $\mu_w$.}
	\end{figure}
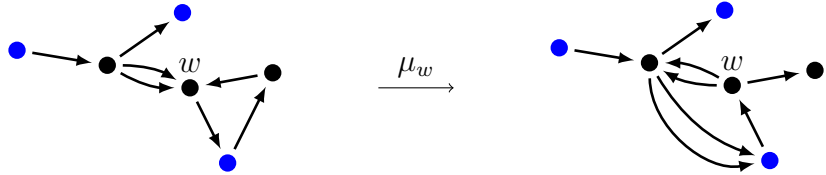
	
	Notice that mutating the same vertex twice in a row gives $Q$ again.
	
	\begin{definition}\label{quiver associated to plg}
		Let $G$ be a plabic graph then we construct a quiver $Q(G)$ as follows. 
		
		\begin{itemize}
			\item[(1)] Place a vertex in any face $f$ of $G$. This vertex is designated as frozen iff $f$ is frozen.
			\item[(2)] For any bicoloured edge $e$ in $G$ separating two faces, we add an arrow between the corresponding vertices directed such that it sees the white endpoint of $e$ to its left.
			\item[(3)] Remove 2-cycles and edges between frozen vertices.
		\end{itemize}
	\end{definition}
	
	We consider the vertices of $Q(G)$ as labelled by the labels of $G$ and write $\mu_I$ for the mutation of $Q(G)$ in direction of a vertex labelled by $I$.
	
	\begin{figure}[!htb]
		\centering
		\begin{tikzpicture}
			\node (A) at ($(-3.5,0)$) {\includegraphics[height = 5cm]{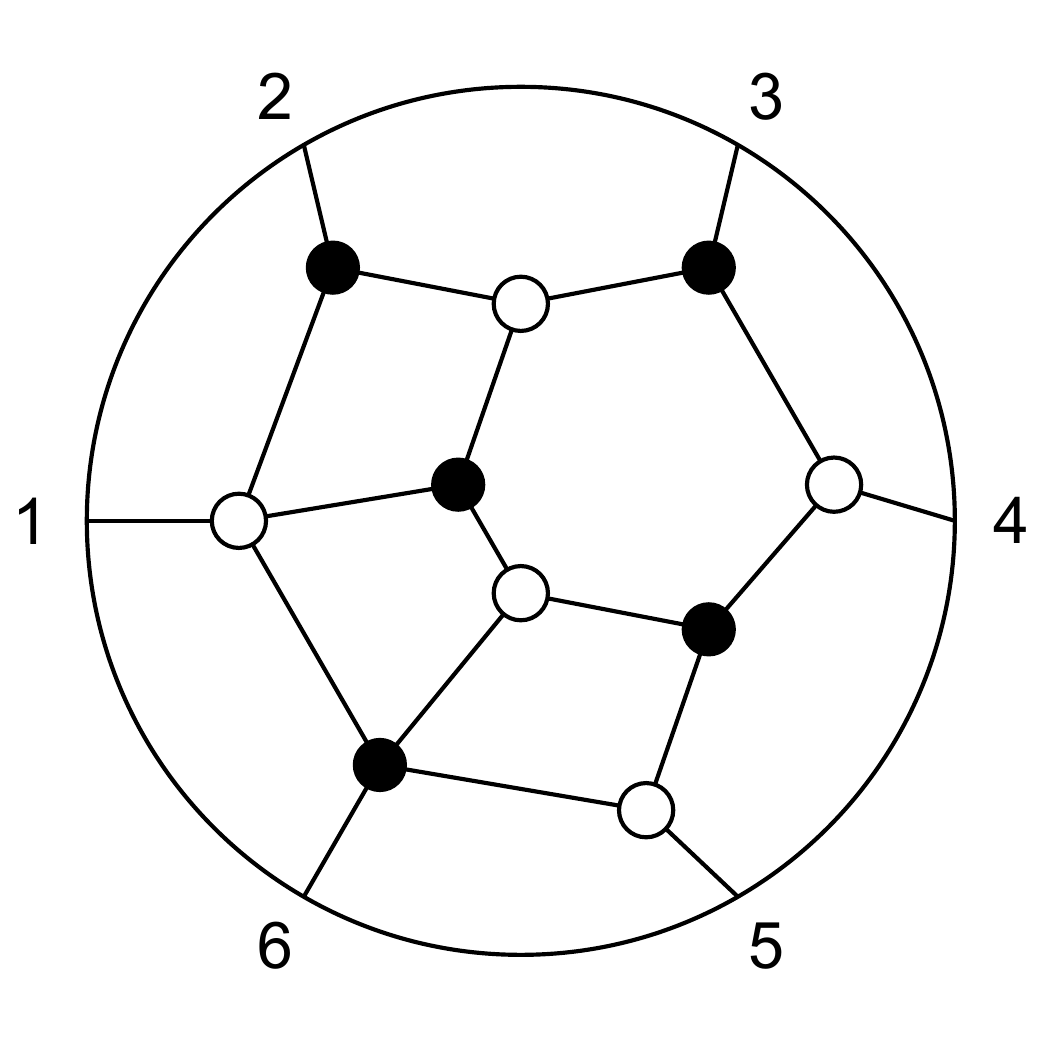}};
			\node (B) at ($( 3.5,0)$) {\includegraphics[height = 5cm]{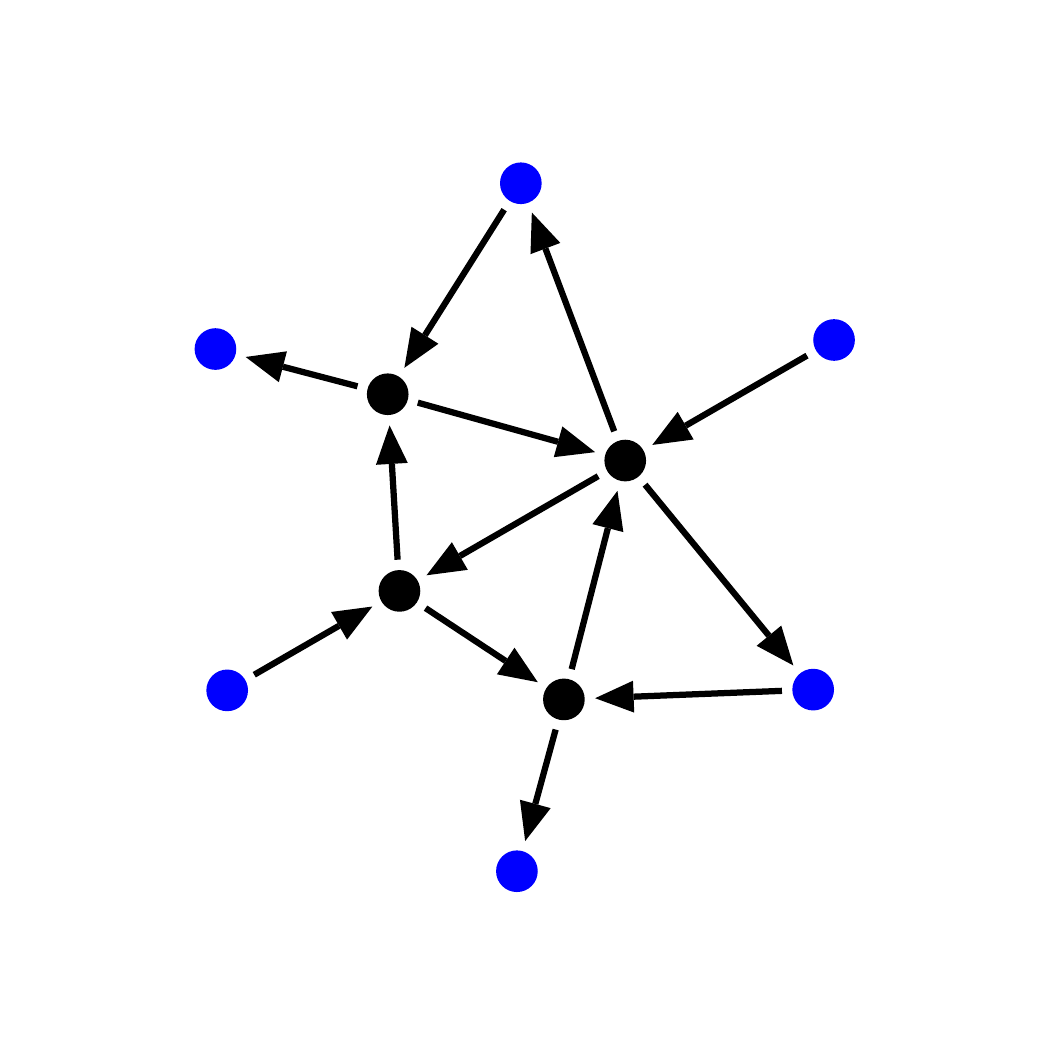}};
			
			\draw[-to, shorten <= 5mm]
			(A) edge (B)
			;
		\end{tikzpicture}
		\caption{\label{fig:figure8} The quiver associated to $G^\text{rec}_{3,6}$.}
	\end{figure}
	
	\begin{lemma}\label{quiver mutation vs plabic mutation}
		The map $G \mapsto Q(G)$ sending a plabic graph to its associated quiver is compatible with the respective notions of mutation.
	\end{lemma}
	
	\section{Grassmannians}
	
	We quickly fix some notation regarding the Grassmannian variety.
	
	\begin{definition}\label{definition grassmannian}
		Let $\text{Gr}_{k,n}$ be the set of $k$-dimensional subspaces of $\C^n$. This set is called the the \textbf{Grassmannian variety}. We will usually use the notation $\mathbb{X}$ for this set where the dependence on $k$ and $n$ is understood.
	\end{definition}
	
	Let $v \in \mathbb{X}$ and $v_1, \ldots, v_k$ be any basis of $v$. If $A_v$ is the matrix whose rows are given by the $v_i$ then $A_v$ is uniquely defined up left multiplication by elements of $\text{GL}_k(\C)$. Thus we identify $\mathbb{X}$ with the set of full rank matrices up to the natural action of $\text{GL}_k(\C)$.
	
	\begin{definition}\label{definition plücker embedding}
		For any $k$-subset $I \subseteq [n]$ we denote the maximal minor with column set $I$ on $M_{k,n}(\mathbb{C})$ by $p_I$.
		The map
		$$ P: \mathbb{X} \to \P^{\binom{n}{k}-1}; v \mapsto ( p_I(A_v) )_I ,$$
		is called the \textbf{Plücker embedding} and the $p_I$ are called \textbf{Plücker coordinates}. 
	\end{definition}
	
	As a shorthand we define $\mathbb{X}^\circ \subset \mathbb{X}$ as the non vanishing set of the coordinates $p_{I_1}, \ldots, p_{I_n}$ corresponding to the frozen labels in any plabic graph. 
	Moreover let $\check{ \mathbb{X}} = \text{Gr}_{n-k, n}$ be the Grassmannian of $n-k$ planes and $\check{ \mathbb{X}}^\circ \subset \check{ \mathbb{X}}$ the non vanishing locus of the Plücker coordinates $p_{J_i}$, where 
	$$J_i = \cyc{i+1}{i+n-k} $$
	is the complement of $I_i$. We usually pass to affine coordinates by normalizing $p_{J_n} = 1$.
	
	\section{$\mathcal{A}$-cluster structure}
	
	In this and in the following sections we will construct face labels of plabic graphs by \enquote{looking to the right} as detailed in \ref{looking to the right}. While cluster algebras were originally introduced by Fomin and Zelevinsky in their foundational paper \cite{Fomin2001}, we once again follow \cite{Rietsch_2019}. \\
	
	Let $G$ be any plabic graph. We define $ \widetilde{ \mathcal{R}}_G$ to be the set of right labels contained in $G$.
	Moreover we set $\mathcal{R}_G = \widetilde{ \mathcal{R}}_G \setminus \{J_n\}.$ \\
	
	With this setup let us have a look at $\mathcal{A}$-seed mutation.
	
	\begin{definition}\label{aseed mutation definition}
		Let $Q$ be a quiver with vertices $V = V(Q)$. We associate a \textbf{cluster variable} $a_v$ to each vertex $v \in V$. If $w \in Q$ is a mutable vertex we define a new set of variables $\mu_w^\mathcal{A}(\{a_v\}) = \{a'_v\}$
		where $a'_v = a_v$ if $v \neq w$, and otherwise, $a'_w$ is determined by
		\begin{equation}\label{aseed mutation equation}
			a_w a'_w = \prod_{v \to w} a_v + \prod_{w \to v} a_v.
		\end{equation}
		Here the products run over all arrows from $v$ to $w$ and $w$ to $v$ in $Q$ respectively. We say that the tuple $(\mu_w(Q), \{a'_v\})$ is obtained from $(Q, \{a_v\})$ by $\mathcal{A}$-\textbf{seed mutation} in direction $w$ and refer to both as (labeled) $\mathcal{A}$-\textbf{seeds}.
	\end{definition}
	
	Note that $\mathcal{A}$-seed mutation is a involution in the sense that mutating in the same direction twice, is the same as doing nothing. We want to use plabic graphs to define seeds, hence we need the following lemma.
	
	\newsavebox{\MIvsLabelsAA}
	\sbox{\MIvsLabelsAA}{
		\resizebox{3cm}{!}{
			\begin{tikzpicture}
				[-,>=stealth',auto, thick]
				\tikzset{white/.style={circle, draw=black, fill=white, scale=1.3}}
				\tikzset{black/.style={circle, draw=black, fill=black, scale=1.3}}
				
				\node[white] (1) at (-0.7, 0.7){};
				\node[black] (2) at ( 0.7, 0.7){};
				\node[white] (3) at ( 0.7,-0.7){};
				\node[black] (4) at (-0.7,-0.7){};
				
				\node (Jac) at ( 0  , 0  ) {$Jac$};
				\node (Jad) at (-1.4, 0  ) {$Jad$};
				\node (Jab) at ( 0  , 1.4) {$Jab$};
				\node (Jbc) at ( 1.4, 0  ) {$Jbc$};
				\node (Jcd) at ( 0  ,-1.4) {$Jcd$};
				
				\coordinate (b1) at ($(0,0)+(135:2.3 and 2.3)$);
				\coordinate (b2) at ($(0,0)+(45:2.3 and 2.3)$);
				\coordinate (b3) at ($(0,0)+(-45:2.3 and 2.3)$);
				\coordinate (b4) at ($(0,0)+(-135:2.3 and 2.3)$);
				
				\draw
				(1) edge (2)
				(2) edge (3)
				(3) edge (4)
				(4) edge (1)
				
				(1) edge (b1)
				(2) edge (b2)
				(3) edge (b3)
				(4) edge (b4)
				;
	\end{tikzpicture} } }
	
	\newsavebox{\MIvsLabelsBB}
	\sbox{\MIvsLabelsBB}{
		\resizebox{3cm}{!}{
			\begin{tikzpicture}
				[-,>=stealth',auto, thick]
				\tikzset{white/.style={circle, draw=black, fill=white, scale=1.3}}
				\tikzset{black/.style={circle, draw=black, fill=black, scale=1.3}}
				
				\node[black] (1) at (-0.7, 0.7){};
				\node[white] (2) at ( 0.7, 0.7){};
				\node[black] (3) at ( 0.7,-0.7){};
				\node[white] (4) at (-0.7,-0.7){};
				
				\node (Jbd) at ( 0  , 0  ) {$Jbd$};
				\node (Jad) at (-1.4, 0  ) {$Jad$};
				\node (Jab) at ( 0  , 1.4) {$Jab$};
				\node (Jbc) at ( 1.4, 0  ) {$Jbc$};
				\node (Jcd) at ( 0  ,-1.4) {$Jcd$};
				
				\coordinate (b1) at ($(0,0)+(135:2.3 and 2.3)$);
				\coordinate (b2) at ($(0,0)+(45:2.3 and 2.3)$);
				\coordinate (b3) at ($(0,0)+(-45:2.3 and 2.3)$);
				\coordinate (b4) at ($(0,0)+(-135:2.3 and 2.3)$);
				
				\draw
				(1) edge (2)
				(2) edge (3)
				(3) edge (4)
				(4) edge (1)
				
				(1) edge (b1)
				(2) edge (b2)
				(3) edge (b3)
				(4) edge (b4)
				;
	\end{tikzpicture} } }
	
	\begin{lemma}\label{three term relation}
		Consider a mutation of a square face in a plabic graph $G$ as shown below.
		\begin{center}
			\begin{tikzpicture}
				\node 			(A) {\usebox{\MIvsLabelsAA}};
				\node at (5,0)	(B) {\usebox{\MIvsLabelsBB}};
				
				\draw[to-to, shorten <=2mm, shorten >=2mm]
				(A) edge (B)
				;
			\end{tikzpicture}
		\end{center}
		Then the corresponding Plücker coordinates of $\check{ \mathbb{X}}$ satisfy the three-term relation
		$$ p_{Jac}p_{Jbd} = p_{Jab}p_{Jcd} + p_{Jbc}p_{Jad}.$$
	\end{lemma}
	
	Let us fix some plabic graph $G$ and define
	$$\mathcal{A}\widetilde{\text{Coord}}_{\check{ \mathbb{X}}}(G) = \{p_J ~|~ J \in \widetilde{ \mathcal{R}}_G\}.$$
	Using the above lemma combined with Theorem \ref{plg_and_wsc_takeaway} we see that any Plücker coordinate of $\check{ \mathbb{X}}$ can be expressed as a rational expression in $\mathcal{A}\widetilde{\text{Coord}}_{\check{ \mathbb{X}}}(G)$. 
	By Theorem \ref{number of faces in a plabic graph} the transcendence degree of $\text{Quot}(\C[\check{ \mathbb{X}}] )$ is precisely $$\dim \check{\mathbb{X}} + 1 = k(n-k) + 1 = |\widetilde{ \mathcal{R}}_G|.$$ 
	Thus the Plücker coordinates in $\mathcal{A}\widetilde{\text{Coord}}_{\check{ \mathbb{X}}}(G)$ 
	must be algebraically independent. 
	Using the language of Definition \ref{aseed mutation definition} we define the $\mathcal{A}$-seed corresponding to $G$ as
	$(Q(G), \mathcal{A}\widetilde{\text{Coord}}_{\check{ \mathbb{X}}}(G) )$. 
	As an easy consequence of Lemma \ref{three term relation} and Theorem \ref{plg_and_wsc_takeaway} we obtain
	
	\begin{proposition}\label{amutation }
		Let $G$ be a plabic graph, $J \in G$ a mutable label and $G' = \mu_J(G)$ then
		$$ \mu_J^\mathcal{A}( \mathcal{A}\widetilde{\text{Coord}}_{\check{ \mathbb{X}}}(G) ) = \mathcal{A}\widetilde{\text{Coord}}_{\check{ \mathbb{X}}}(G') .$$
	\end{proposition}
	
	As consequence of the well known Laurent phenomenon (see for example \cite{Fomin2001} or \cite{Fomin2016}) and its positivity (see \cite{Lee2013}) we have the following important fact.
	
	\begin{proposition}\label{laurent phenomenon}
		Any Plücker coordinate of $\check{ \mathbb{X}}$ can be expressed as a Laurent polynomial in terms of the cluster variables in any $\mathcal{A}$-seed which is mutation equivalent to an $\mathcal{A}$-seed coming from a plabic graph. 
		Moreover all appearing coefficients are positive integers.
	\end{proposition}
		
	\section{The Superpotential and Newton-Okounkov bodies}
	In their paper \cite{Rietsch_2019}, Rietsch and Williams associate Newton-Okounkov bodies to plabic graphs and show how they can be obtained by tropicalizing a certain rational function. We review this approach here.
	
	\begin{definition}
		Define $J_i^+ = \cyc{i+1}{i+n-k-1} \cup \{i+n-k+1\}$ and
		let $W_i = p_{J_i}^+/p_{J_i}$. We call the regular function $W: \check{ \mathbb{X}}^\circ \times \C^* \to \C$ defined by
		$$ W = \sum_{i = 1}^n W_i q^{\delta_{i,k}} ,$$
		the \textbf{superpotential} $W$ on $\check{ \mathbb{X}}$. Here $\delta_{i,k}$ denotes the Kronecker delta, and $q$ is the coordinate $\C^*$ in the product above product of varieties.
	\end{definition}
	
	By Proposition \ref{laurent phenomenon} and Theorem \ref{plg_and_wsc_takeaway} we can always rewrite the superpotential as a Laurent polynomial with positive integer coefficients in the variables of a given plabic graph.
	For rectangle graphs a closed formula is known (see Proposition 10.5 of \cite{Rietsch_2019}).
	
	\begin{definition}\label{tropicalization}
		Let $h = \sum_{u \in \Z^m} c_u X^u$ be a Laurent polynomial in the variables $X_1, ..., X_m$ with non-negative coefficients. We define the \textbf{tropicalization} $\text{Trop}(h): \R^m \to \R$ of $h$ to be the piecewise linear map 
		$$ \text{Trop}(h)(v) = \min( v \cdot u ~|~ c_u \neq 0 ) ,$$
		where $v \cdot u$ denotes the standard scalar product on $\R^m$.
	\end{definition}
	
	Given any plabic graph $G$ we can now apply the tropicalization procedure defined above to obtain a piecewise linear map 
	$$ \text{Trop}_G(W): \R^{\mathcal{R}_G} \times \R \to \R .$$
	
	Here the additional factor $\mathbb{R}$ corresponds to the coordinate $q$ in the superpotential. 
	Since $p_{J_n}$ does not vanish on $\check{ \mathbb{X}}^\circ$ we usually normalize $p_{J_n} = 1$. This explains why we
	obtain a function on $\R^{\mathcal{R}_G} \times \R$ instead of a function on $\R^{\widetilde{\mathcal{R}_G}} \times \R$.
	
	\begin{definition}\label{super potential polytope}
		Let $G$ be any plabic graph and $\text{Trop}_G(W)$ be constructed as explained above. For any $r \in \R$ the polyhedron
		$$\Gamma_G^r = \{v \in \R^{\mathcal{R}_G} ~\mid~ \text{Trop}_G(W)(v,r) \geq 0\} ,$$ 
		is called the \textbf{superpotential polytope}. We abbreviate $\Gamma_G^1$ by $\Gamma_G$.
	\end{definition}
	
	The word \textit{polytope} above is not an accident, although it is not immediate from the definition that we are indeed dealing with polytopes here (see Proposition \cite{Rietsch_2019}). 
	
	\begin{theorem}[\cite{Rietsch_2019} Definition 8.2, Remark 8.5 and Theorem 16.18]
		For any plabic graph $G$ the superpotential polytope $\Gamma_G$ defined above equals the associated Newton-Okounkov body $\Delta_G$.
	\end{theorem}
	
	Let us denote the Gelfand-Tsetlin polytope with top row $(0^k,r^{n-k})$ by $\mathcal{GT}_{k,n}^r$ (see \cite{Gelfand1950} or \cite{Alexandersson2020} for more details).
	
	\begin{proposition}\label{nob is gt for rec}
		Let $G = G^\text{rec}_{k,n}$. Then the Newton-Okounkov body $\Delta_G$ is unimodular equivalent to $\mathcal{GT}_{k,n}^1$.
	\end{proposition}
	
	We further denote by $\textnormal{FFLV}_{k,n}^r$ the Feigin-Fourier-Littelmann-Vinberg polytope (see \cite{FF2018} and \cite{Feigin2011}). 
	
	\begin{theorem}
		Let $G = \widehat G^{\text{rec}}_{k,n}$. Then the Newton-Okounkov body $\Delta_G$ is unimodular equivalent to $\textnormal{FFLV}_{k,n}^1$. 
	\end{theorem}
	
	Since $\widehat G^{\text{rec}}_{k,n}$ is a reflection of $ G^{\text{rec}}_{k,n}$ this means that plabic graphs in the same orbit will in general still admit non equivalent Newton-Okounkov bodies.
	
	\section{The superpotential and (dual) checkboards}
	
	We now derive closed formulas for the superpotential in terms of (dual) checkboard graphs.
	This more or less amounts to finding the labels $J_i^+$ starting from $G^\text{ch}_{k,n}$, and then going through the corresponding seed mutations.
	
	\begin{example}\label{super labels check example}
		Consider $G = G^\text{ch}_{3,6}$. Two of the desired labels, namely $J_1^+$ and $J_4^+$, are already contained in $G$. Indeed 
		$ J_1^+ = J^\text{ch}_{1,2} = \{2,3,5\}$ and $J_4^+ = J^\text{ch}_{2,1} = \{2,5,6\} ,$
		where $J^\text{ch}_{i,j}$ denotes the complement of $I^\text{ch}_{i,j}$ (see Proposition \ref{checkboard_graph_labels}). \\
		\begin{figure}[!htb]
			\centering
				\includegraphics[height=6cm]{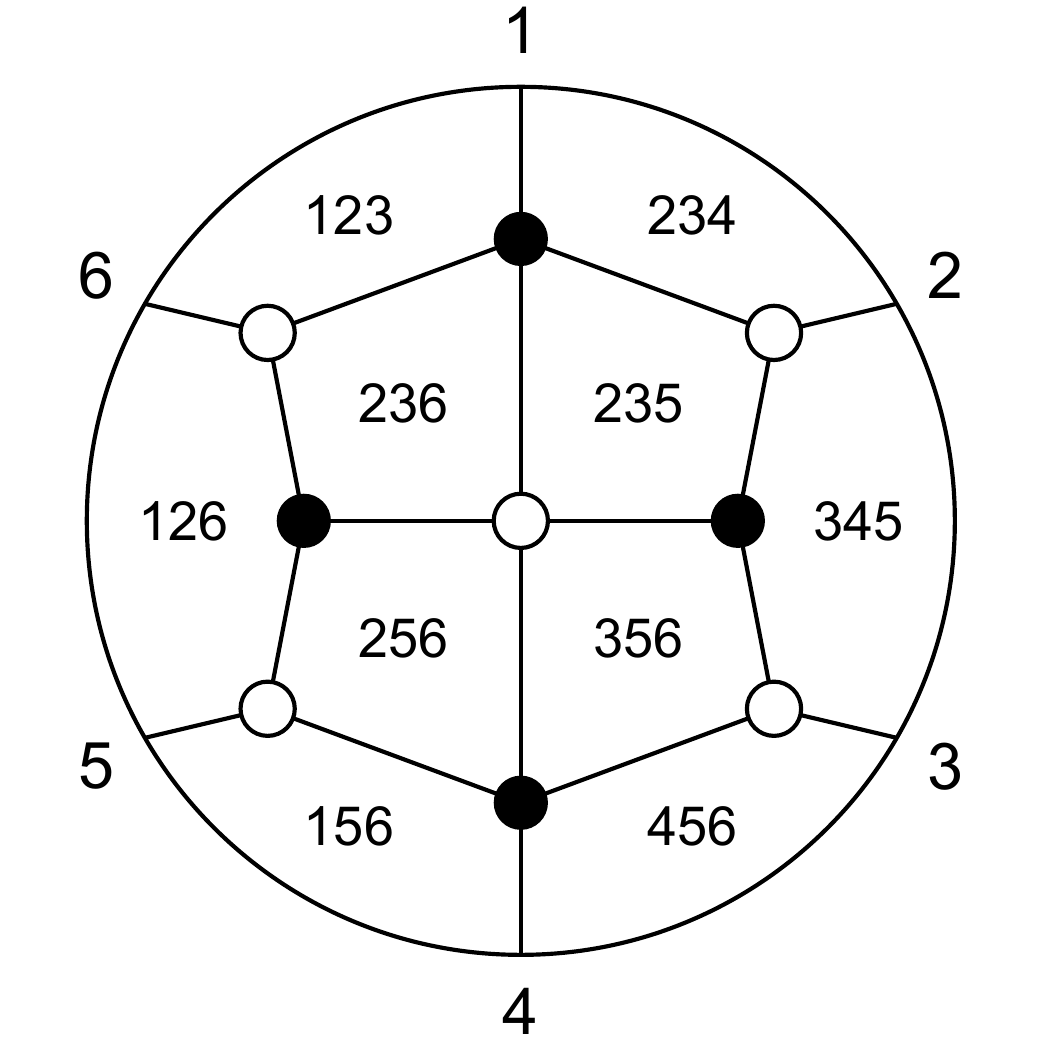}
			\caption{\label{fig:figure9} The checkboard graph $G^\text{ch}_{3,6}$.}
		\end{figure}
		
		One easily verifies that mutating $J^\text{ch}_{i,j}$ with odd $i+j$ will exchange it with $\sigma J^\text{ch}_{i,j}$. Thus we easily find $J_2^+$ and $J_5^+$. Finding $J_3^+$ and $J_6^+$ requires a little more effort.
		Let us keep track of the original grid positions of faces during mutations. That is, we refer to the inner face in row $i$ and column $j$ of $G$ as $f_{i,j}$ even after the grid is deformed by mutations. Then, to obtain $J_3^+$ we mutate $f_{2,1}, f_{1,1}$ and $f_{2,2}$ in that order. Similarly to obtain $J_6^+$ we may mutate the faces $f_{1,2}, f_{2,2}$ and $f_{1,1}$. In each case the face mutated last will carry the desired label. 
	\end{example}
	
	The above example generalizes to arbitrary $k$ and $n$. 
	Let us fix some notation. We call the set of inner faces $f_{i,j}$ in $G^\text{ch}_{k,n}$ whose indices satisfy $i-j = d$, the \textbf{$d$-th diagonal} of $G^\text{ch}_{k,n}$. \\
	
	Now for any $d \in \{2-k, \ldots, n-k-2\}$ such that $d$ is even, we have a \textbf{downward diagonal sequence} of mutations given by first mutating all faces of the $d+1$-th diagonal (if there is such a diagonal) in any order followed by the faces in the $d$-th diagonal from top to bottom. Below the downward diagonal sequence for $d = -2$ in $G^\text{ch}_{5,9}$ is illustrated, where positions sharing the same number are to be mutated in any order.
	$$
	\begin{array}{|c|c|c|c|}
		\hline 
		\rule{0pt}{2.5ex}
		\hspace{1ex} & 1 & 2 &   \\ \hline \rule{0pt}{2.5ex}
		& & 1 & 3  \\ \hline \rule{0pt}{2.5ex}
		& &   & 1 \\ \hline 
	\end{array}
	$$
	Similarly we have an \textbf{upward diagonal sequence} given by first mutating all faces the $d-1$-th diagonal followed by mutating the faces in the $d$-th diagonal from bottom to top. Below the upwards diagonal sequence for $d = 0$ in $G^\text{ch}_{5,9}$ is illustrated.
	$$
	\begin{array}{|c|c|c|c|}
		\hline 
		\rule{0pt}{2.5ex}
		4 & 1 & &   \\ \hline \rule{0pt}{2.5ex}
		& 3 & 1 &   \\ \hline \rule{0pt}{2.5ex}
		& & 2   & 1 \\ \hline 
	\end{array}
	$$
	
	\begin{lemma}\label{diagonal sequences}
		Let $G = G^\text{ch}_{k,n}$, then we always have
		$$
		J^\text{ch}_{1,k-1} = J_{\lfloor k/2 \rfloor}^+ 
		~~~\text{and}~~~
		J^\text{ch}_{n-k-1,1} = J_{\lfloor (n+k)/2 \rfloor}^+ . 
		$$
		Moreover if $d \in \{2-k, \ldots, n-k-2\}$ is even, then after applying the corresponding downward diagonal sequence, the last face mutated will carry the label $\sigma^{d/2} J_k^+$. 
		Applying the upward diagonal sequence instead gives the label $\sigma^{-d/2} J_n^+$. \\
		Finally if $k$ or $n-k$ is odd then mutating $f_{1,k-1}$ or $f_{n-k-1,1}$ will yield the labels $J_{\lceil k/2 \rceil}^+ $ or $J_{\lceil (n+k)/2 \rceil}^+$ respectively.
	\end{lemma}
	
	\begin{proof}
		The first claim follows directly from Proposition \ref{checkboard_graph_labels} by taking complements. Explicitly we have 
		$$ J_{i,j} = J^\text{ch}_{i,j} = \sigma^{- \lceil (i+j)/2 \rceil}( [i] \cup [n-k+j] \setminus [i+j] ) ,$$ where we also allow $i \in \{0, n-k\}$ and $j \in \{0,k\}$.
		If $k$ is odd, mutating $f_{1,k-1}$ will rotate its label $J_{\lfloor k/2 \rfloor}^+$ clockwise i.e. exchange it by $\sigma J_{\lfloor k/2 \rfloor}^+ = J_{\lceil k/2 \rceil}^+$. Similarly one can verify the special case for odd $n-k$. The diagonal sequences are dealt with inductively. We omit the details here for the sake of brevity.
	\end{proof}
	
	It is not too hard to see that the above lemma exhausts all possible $J_i^+$. 
	The next step is to consider the $\mathcal{A}$-seed mutations corresponding to the above sequences. 
	For this purpose set
	$$ D = \left\{ - \left\lfloor \frac{k-2}{2} \right\rfloor, \ldots,  \left\lfloor \frac{n-k-2}{2} \right\rfloor  \right\} .$$
	Moreover define
	$y_{i,j} = y_{i,j}^\text{ch}$ to be the Plücker coordinate of $\check{ \mathbb{X}}$ indexed by $J^\text{ch}_{i,j}$.
	Note that some of the $y_{i,j}$ coincide. For example $y_{0,0} = y_{0,1}$ (see Remark \ref{extended checkboard grid}).
	
	\begin{lemma}\label{superpotetial terms for check}
		For any 
		$c \in \{- \left\lceil \frac{k-2}{2} \right\rceil \} \cup D $
		we have the equality
		\begin{equation}\label{downward diag formula}
			\frac{ p_{J_{ k+c }^+ } }{ p_{J_{ k+c } } } = 
			\sum_{a-b = 2c} \frac{ y_{a,b+1}y_{a+1,b-1} }{ y_{a,b}y_{a+1,b} } +
			\sum_{a-b = 2c} \frac{ y_{a-1,b}y_{a+1,b-1} }{ y_{a,b-1}y_{a,b} },
		\end{equation}
		where the $a,b$ in the above sums run over all possibilities that make the occurring variables defined. 
		Similarly, for any 
		$c \in D \cup \{ \left\lceil \frac{n-k-2}{2} \right\rceil \}$
		we have
		\begin{equation}\label{upward diag formula}
			\frac{ p_{J_{ n-c }^+ } }{ p_{J_{ n-c } } } = 
			\sum_{a-b = 2c-1} \frac{ y_{a,b+1}y_{a+1,b-1} }{ y_{a,b}y_{a+1,b} } +
			\sum_{a-b = 2c-1} \frac{ y_{a-1,b}y_{a+1,b-1} }{ y_{a,b-1}y_{a,b} }.
		\end{equation}
		Both equalities should be viewed as equalities of regular functions on the open subset of $ \check{\mathbb{X}} $ where the $y_{i,j}$ do not vanish.
	\end{lemma}
	
	Putting everything together, we obtain a formula of the superpotential $W$ in terms of the coordinates of $G^\text{ch}_{k,n}$.
	
	\begin{proposition}\label{W in terms of check formula}
		Let $G = G^\text{ch}_{k,n}$ be the checkboard graph. Then on the open subset of $\check{ \mathbb{X}}^\circ$ where all $y_{i,j} \neq 0$, the superpotential $W$ equals
		\begin{equation*}
			W = \frac{y_{1,k-1}}{y_{0,k}} + \frac{y_{n-k-1,1}}{y_{n-k,0}} 
			+ \sum_{d = 1-k}^{n-k-2} 
			\left(~ 
			\sum_{a-b = d} \frac{y_{a,b+1} ~ y_{a+1,b-1}}{y_{a,b} ~ y_{a+1,b}} +
			\sum_{a-b = d} \frac{y_{a-1,b} ~ y_{a+1,b-1}}{y_{a,b-1} ~ y_{a,b}}
			\right)q^{\delta_{0,d}}
		\end{equation*} 
		as regular function. Here $y_{i,j}$ denotes the Plücker coordinate of $\check{ \mathbb{X} }$ whose label is given by $$ J^\text{ch}_{i,j} = \sigma^{- \lceil (i+j)/2 \rceil}( [i] \cup [n-k+j] \setminus [i+j] ) .$$
	\end{proposition}
	
	Notice that this formula looks exceptionally similar to the one for the rectangle graphs. One should keep in mind though that the $y_{i,j}$ here satisfy different equalities as for rectangle graph (compare Remark \ref{extended rectangle grid} and Remark \ref{extended checkboard grid}).
	
	What about the graphs in the orbit of $G^\text{ch}_{k,n}$? As a first step note that $\sigma^m J_i = J_{i+m}$ and $\sigma J_i^+ = J_{i+m}^+$ for any $m$. Then combining Lemma \ref{mutation_vs_action} with the lemmata \ref{diagonal sequences} and \ref{superpotetial terms for check} we obtain the following proposition.
	
	\begin{proposition}\label{W in terms of rotated check formula}
		Let $m \in [n]$ and $G = \sigma^m G^\text{ch}_{k,n}$ be a rotation of the checkboard graph. Moreover let
		$y_{i,j}$ denote the Plücker coordinate of $\check{ \mathbb{X} }$ labelled by $\sigma^m J^\text{ch}_{i,j}$.
		Then on the open subset of $\check{ \mathbb{X}}^\circ$ where all $y_{i,j} \neq 0$, the superpotential $W$ equals
		\begin{align*}
			W = &~ \frac{y_{1,k-1}}{y_{0,k}}q^{ \delta_{m, \left\lceil k/2 \right\rceil } } 
			+ \frac{y_{n-k-1,1}}{y_{n-k,0}}q^{ \delta_{m, \left\lceil (n+k)/2 \right\rceil } } \\
			& + \sum_{d = 1-k}^{n-k-2}
			\left(~
			\sum_{a-b = d} \frac{y_{a,b+1} ~ y_{a+1,b-1}}{y_{a,b} ~ y_{a+1,b}} +
			\sum_{a-b = d} \frac{y_{a-1,b} ~ y_{a+1,b-1}}{y_{a,b-1} ~ y_{a,b}}
			\right)
			q^{\epsilon_{m,d}}
		\end{align*}
		as a regular function. Here 
		$$ \epsilon_{m, d} = \delta_{m, -d/2} + \delta_{m, (d+1)/2 + k} ,$$ 
		and of course $\sigma = x \mapsto x+1$.
	\end{proposition}
	
	We add that the arguments of the Kronecker deltas above should be interpreted modulo $n$, wherever possible.
	
	\begin{proof}
		We obtain the same formulas for the terms of the superpotential as we did when considering $G^\text{ch}_{k,n}$. They are just shifted cyclically. I.e. the formula for say $p_{J_i^+}/p_{J_i}$ in terms of $G^\text{ch}_{k,n}$ will now provide the formula for $p_{J_{i+m}^+}/p_{J_{i+m}}$. Thus we only need to correct the position of the extra variable $q$ accordingly.
	\end{proof}
	
	Since the labels $J_i^+$ are not invariant under reflections, we cannot use the same idea to obtain a formulas for $W$ in the case where $G$ is a reflection $G^\text{ch}_{k,n}$. One can however repeat the same analysis as for $G^\text{ch}_{k,n}$ in order to obtain formulas anyway.
	
	\begin{proposition}\label{W in terms of reflected check formula}
		Let $m \in [n]$ and $G = \sigma^m \rho G^\text{ch}_{k,n}$ be a reflection of the checkboard graph. Moreover let
		$y_{i,j}$ denote the Plücker coordinate of $\check{ \mathbb{X} }$ labelled by
		$$ 
		\sigma^m \rho J^\text{ch}_{j, i} = 
		\sigma^{ \left\lceil \frac{i+j}{2} \right\rceil + m}( [n] \setminus [n-j] \cup [n-j-i] \setminus [k-i] ) ,
		$$
		where $i \in \{0, \ldots, k\}$ and $j \in \{0, \ldots, n-k\}$.
		Then on the open subset of $\check{ \mathbb{X}}^\circ$ where all $y_{i,j} \neq 0$, the superpotential $W$ equals
		\begin{align*}
			W = &~ \frac{y_{1,k-1}}{y_{0,k}}q^{ \delta_{m, 0 } } 
			+ \frac{y_{n-k-1,1}}{y_{n-k,0}}q^{ \delta_{m, \left\lfloor n/2 \right\rfloor } } \\
			& + \sum_{d = 1}^{n-2}
			\left(~ 
			\sum_{a+b = d} \frac{y_{a,b-1} ~ y_{a+1,b+1}}{y_{a,b} ~ y_{a+1,b}} +
			\sum_{a+b = d} \frac{y_{a-1,b} ~ y_{a+1,b+1}}{y_{a,b} ~ y_{a,b+1}}
			\right)
			q^{\epsilon_{m,d}}
		\end{align*} 
		as regular function. Here 
		$$ \epsilon_{m, d} = \delta_{m, d/2} + \delta_{m, -(d+1)/2} ,$$ 
		and $\sigma = x \mapsto x+1$, $\rho = x \mapsto n+1-x$.
	\end{proposition}
	
	An important difference to the case of $G^\text{ch}_{k,n}$ is that the diagonals used for diagonal sequences are defined as the set of faces $f_{i,j}$ such that $i+j = d$ for a fixed $d$. Intuitively this means these diagonals stretch from the top right to the bottom left. \\
	
	Essentially the same analysis can also be applied to dual checkboard graphs $\widehat G^\text{ch}_{k,n}$.
	
	\begin{proposition}\label{W in terms of rotated dual check formula}
		Let $m \in [n]$ and $G = \sigma^m \widehat G^\text{ch}_{k,n}$ be a rotation of the dual checkboard graph. Moreover let
		$y_{i,j}$ denote the Plücker coordinate of $\check{ \mathbb{X} }$ labelled by
		$$ 
		\sigma^m \hat J^\text{ch}_{i, j} = 
		\sigma^{ - \left\lceil \frac{i+j}{2} \right\rceil + m}( [n] \setminus [k+j] \cup [i+j] \setminus [i] ) ,
		$$
		where $i \in \{0, \ldots, k\}$ and $j \in \{0, \ldots, n-k\}$.
		Then on the open subset of $\check{ \mathbb{X}}^\circ$ where all $y_{i,j} \neq 0$, the superpotential $W$ equals
		\begin{align*}
			W = &~ \frac{y_{1,k-1}}{y_{0,k}}q^{ \delta_{m, 0 } } 
			+ \frac{y_{n-k-1,1}}{y_{n-k,0}}q^{ \delta_{m, \left\lceil n/2 \right\rceil } } \\
			& + \sum_{d = 1}^{n-2}
			\left(~
			\sum_{a+b = d} \frac{y_{a,b-1} ~ y_{a+1,b+1}}{y_{a,b} ~ y_{a+1,b}} +
			\sum_{a+b = d} \frac{y_{a-1,b} ~ y_{a+1,b+1}}{y_{a,b} ~ y_{a,b+1}}
			\right)
			q^{\epsilon_{m,d}}
		\end{align*} 
		as regular function. Here 
		$$ \epsilon_{m, d} = \delta_{m, -d/2} + \delta_{m, (d+1)/2} ,$$ 
		and $\sigma = x \mapsto x+1$.
		
	\end{proposition}
	
	Similar for reflections of $\widehat G^\text{ch}_{k,n}$, we obtain
	
	\begin{proposition}\label{W in terms of reflected dual check formula}
		Let $m \in [n]$ and $G = \sigma^m \rho \widehat G^\text{ch}_{k,n}$ be a reflection of the dual checkboard graph. Moreover let
		$y_{i,j}$ denote the Plücker coordinate of $\check{ \mathbb{X} }$ labelled by
		$$ 
		\sigma^m \rho \hat J^\text{ch}_{j, i} = 
		\sigma^{ \left\lceil \frac{i+j}{2} \right\rceil + m }( [n-k-i] \cup [n-j] \setminus [n-i-j] ) ,
		$$
		where $i \in \{0, \ldots, n-k\}$ and $j \in \{0, \ldots, k\}$.
		Then on the open subset of $\check{ \mathbb{X}}^\circ$ where all $y_{i,j} \neq 0$, the superpotential $W$ equals
		\begin{align*}
			W = &~ \frac{y_{1,k-1}}{y_{0,k}}q^{ \delta_{m, \left\lfloor k/2 \right\rfloor } } 
			+ \frac{y_{n-k-1,1}}{y_{n-k,0}}q^{ \delta_{m, \left\lfloor (n+k)/2 \right\rfloor } } \\
			& + \sum_{d = 1-k}^{n-k-2}
			\left(~
			\sum_{a-b = d} \frac{y_{a,b+1} ~ y_{a+1,b-1}}{y_{a,b} ~ y_{a+1,b}} +
			\sum_{a-b = d} \frac{y_{a-1,b} ~ y_{a+1,b-1}}{y_{a,b-1} ~ y_{a,b}}
			\right)
			q^{\epsilon_{m,d}}.
		\end{align*}
		as a regular function. Here 
		$$ \epsilon_{m, d} = \delta_{m, d/2+k} + \delta_{m, -(d+1)/2} ,$$ 
		and $\sigma = x \mapsto x+1$, $\rho = x \mapsto n+1-x$.
		
	\end{proposition}
	
	\section{Newton-Okounkov bodies of (dual) checkboards}
	
	We are now ready to investigate Newton-Okounkov bodies.
	
	\begin{remark}\label{inequalities for rotated check}
		Let $G = \sigma^m G^\text{ch}_{k,n}$. Then the defining inequalities of $\Delta_G \subseteq \R^{\mathcal{R}_G} \cong \R^{k(n-k)}$ are
		\begin{align}
			\delta_{m, \left\lceil k/2 \right\rceil } 	  & \geq v_{0,k} - v_{1,k-1}, \nonumber \\ 
			\delta_{m, \left\lceil (n+k)/2 \right\rceil } & \geq v_{n-k,0} - v_{n-k-1,1}, \nonumber \\ 
			\epsilon_{m, a-b} & \geq v_{a,b} - v_{a+1,b-1} - v_{a,b+1} + v_{a+1,b}, \\
			\epsilon_{m, a-b} & \geq v_{a,b} - v_{a+1,b-1} - v_{a-1,b} + v_{a,b-1}, 
		\end{align}
		where $$ \epsilon_{m, d} = \delta_{m, -d/2} + \delta_{m, (d+1)/2 + k} $$
		and we have an inequality of the form (4) or (5) for any $a \in \{0, \ldots, n-k-1 \}, b \in \{1, \ldots, k\}$ that make the involved variables defined.
	\end{remark}
	
	\begin{proposition}\label{nobs are gt for rotated check}
		Let $G = \sigma^m G^\text{ch}_{k,n}$, then the Newton-Okounkov body $\Delta_G \subseteq \R^{\mathcal{R}_G} \cong \R^{k(n-k)}$ of $G$ is unimodular equivalent to the Gelfand-Tsetlin polytope $\mathcal{GT}_{k,n}^1$, with defining inequalities given by
		\begin{align*}
			1 & \geq f_{0,k}, \\ 
			0 & \geq -f_{n-k-1,1}, \\ 
			0 & \geq f_{a,b} - f_{a,b+1}, \\
			0 & \geq f_{a,b} - f_{a-1,b}.
		\end{align*}
	\end{proposition}
	
	\begin{proof}
		The idea is a variation of \cite{Rietsch_2019} lemma 16.2: We want to define a unimodular transformation by setting $f_{i,j} = v_{i,j} - v_{i+1, j-1}$, such that the inverse is given by $f_{i,j} + f_{i+1,j-1} + \ldots $. 
		However, if we try this then we run into trouble as $f_{i,j} + f_{i+1,j-1} + \ldots = v_{i,j} + v_{a,b}$ with $a = n-k$ or $b = 0$ and $v_{a,b} = 0$ is usually wrong.
		
		To circumvent this problem we look at the polytope in a suitable subspace of $\R^{(k+1)(n-k+1)}$. Recall that the $v_{i,j}$ satisfy some equalities like $v_{0,0} = v_{0,1}$, stemming from the boundary faces of $G$ (see also Remark \ref{extended checkboard grid}). So far this was just a consequence of our notation, i.e. $v_{0,0}$ and $v_{0,1}$ were never separate variables to begin with.
		
		We will separate them now. For this purpose let $w_{i,j}$ be the coordinates of $\R^{(k+1)(n-k+1)}$. 
		We define a subspace $U_1$ of $\R^{(k+1)(n-k+1)}$ by making the equalities between the $v_{i,j}$ explicit. That is, we require
		\begin{equation*}
			\begin{array}{cc}
				w_{i,0} = w_{i-1,0} 	& (2 ~|~ i) \\
				w_{n-k,j} = w_{n-k,j-1} & (2 ~|~ n-k+j) \\
				w_{0,j} = w_{0,j+1} 	& (2 ~|~ j) \\
				w_{i,k} = w_{i+1,k} 	& (2 ~|~ k+i),
			\end{array}
		\end{equation*}
		where additionally, depending on $m$, the variables in exactly one of the above equations should be set to zero, since we normalized $p_{J_n} = 1$.
		More precisely, we have $v_{i,j} = 0$ if the corresponding Plücker coordinate $y_{i,j}$ (defined as in Proposition \ref{W in terms of rotated check formula}) carries the label $J_n$, and in this case we also want to have $w_{i,j} = 0$. 
		
		It will not matter for our proof which of the above $w_{i,j}$ exactly are set to zero, just that two of them are, and that they correspond to the same equation or boundary face of $G$.  
		
		Note that we have an obvious inclusion $\iota_1: \R^{k(n-k)} \to \R^{(k+1)(n-k+1)}$ that takes $\Delta_G$ to the polytope inside $U_1$ defined by the inequalities 
		\begin{align*}
			\delta_{m, \left\lceil k/2 \right\rceil } 	  & \geq w_{0,k} - w_{1,k-1}  \\ 
			\delta_{m, \left\lceil (n+k)/2 \right\rceil } & \geq w_{n-k,0} - w_{n-k-1,1}  \\ 
			\epsilon_{m, a-b} & \geq w_{a,b} - w_{a+1,b-1} - w_{a,b+1} + w_{a+1,b} \\
			\epsilon_{m, a-b} & \geq w_{a,b} - w_{a+1,b-1} - w_{a-1,b} + v_{a,b-1}
		\end{align*}
		Of course $\iota_1$ is an isomorphism onto $U_1$, whose inverse is obtained by restricting the appropriate projection $\pi_1: \R^{(k+1)(n-k+1)} \to \R^{k(n-k)}$. Note that both $\iota_1$ and $\pi_1$ send lattice points onto lattice points. 
		Next, we define a linear map $\tilde F:\R^{(k+1)(n-k+1)} \to \R^{(k+1)(n-k+1)}$ by sending $w_{i,j}$ to 
		$$
		f_{i,j} = \begin{cases}
			w_{i,j} - w_{i+1,j-1} & \text{if}~ i < n-k ~\text{and}~ j > 0 \\
			w_{i,j}				  & \text{otherwise}.
		\end{cases}
		$$
		The second case in the above definition fixes our earlier issue, as we now have 
		$$ f_{i,j} + f_{i+1,j-1} + \ldots = w_{i,j} $$ again. Thus $\tilde F$ is a unimodular transformation. \\
		
		Rewriting the (in)equalities in $f$-coordinates gives 
		\begin{align*}
			\delta_{m, \left\lceil k/2 \right\rceil } 	  	& \geq f_{0,k}  			\\ 
			\delta_{m, \left\lceil (n+k)/2 \right\rceil } 	& \geq -f_{n-k-1,1}  		\\ 
			\epsilon_{m, a-b} 								& \geq f_{a,b} - f_{a,b+1} 	\\
			\epsilon_{m, a-b} 								& \geq f_{a,b} - f_{a-1,b}, 
		\end{align*}
		together with
		\begin{equation}\label{nobs are gt equations 1}
			\begin{array}{cccc}
				f_{i,0} 						&=& f_{i-1,0} 							& (2 ~|~ i)     \\
				f_{n-k,j} 						&=& f_{n-k,j-1} 						& (2 ~|~ n-k+j) \\
				f_{0,j} + f_{1,j-1} + \ldots  	&=& f_{0,j+1} + f_{1,j} + \ldots 		& (2 ~|~ j)     \\
				f_{i,k} + f_{i+1,k-1} + \ldots 	&=& f_{i+1,k} + f_{i+2,k-1} + \ldots 	& (2 ~|~ k+i)
			\end{array}
		\end{equation}
		where again one of the above equalities should additionally be equal to zero. We can rewrite the latter as 
		\begin{equation}\label{nobs are gt equations 2}
			\begin{array}{cccc}
				f_{i,0} - f_{i-1,0} 			&=& 0 	& (2 ~|~ i)       	\\
				f_{n-k,j} - f_{n-k,j-1} 		&=& 0 	& (2 ~|~ n-k+j)   	\\
				f_{i,0} - f_{i-1,0} + s_i 		&=& 0 	& (2 ~\nmid~ i)   	\\
				f_{n-k,j} - f_{n-k,j-1} + t_j 	&=& 0 	& (2 ~\nmid~ n-k+j),
			\end{array}
		\end{equation}
		where $s_i$ and $t_j$ are linear combinations of the $f_{i,j}$ with $i < n-k$ and $j > 0$. Let us now fix an order of the coordinates $f_{i,j}$ such that $$f_{0,0} < f_{1,0} < \ldots < f_{n-k,0} < f_{n-k, 1} < \ldots < f_{n-k,k}$$ and $f_{i,j} < f_{0,0}$ for any $i < n-k$ and $j >0$.
		Taking into account the zero equality, we can represent the equations (\ref{nobs are gt equations 2}) 
		in the form $B f = 0$ where 
		$$B = [B_1 | B_2] \in \Z^{ n+1 \times (k+1)(n-k+1) } ,$$ such that $B_2 \in \Z^{n+1 \times n+1}$. 
		More precisely, let $e_i$ denote the $i$-th unit (row-)vector of $\R^{n+1}$. Possibly after a reordering of equalities, the rows of $B_2$ are given by
		$$ e_2-e_1, ~e_3-e_2, ~\ldots~ , ~e_l-e_{l-1}, ~e_l, ~e_{l+1}, ~e_{l+2}-e_{l+1}, ~\ldots~ , ~e_{n+1}-e_n ,$$
		for some $l$ depending on which of the equalities (\ref{nobs are gt equations 1}) is furthermore equal to zero. Important for us is that in any case we have $B_2 \in \text{GL}_{n+1}(\Z)$. \\
		Now let $\pi_2: \R^{(k+1)(n-k+1)} \to \R^{k(n-k)}$ be the projection that kills any $f_{i,j}$ with $i = n-k$ or $j = 0$. We claim that $\pi_2$ restricts to an isomorphism $\pi_2:U_2 \to \R^{k(n-k)}$, 
		where $U_2 \subset \R^{(k+1)(n-k+1)}$ is the kernel of $B$.
		Indeed the inverse is given by
		$$\iota_2: \R^{k(n-k)} \to U_2; f \mapsto (f, -B_2^{-1}B_1f).$$ 
		Again, both $\iota_2$ and $\pi_2$ send lattice points onto lattice points. As a consequence the composition 
		$$ F = \pi_2 \tilde F \iota_1 $$ 
		is a unimodular transformation. By construction $F$ takes $\Delta_G$ to the polytope in $\R^{k(n-k)}$ defined just by the inequalities 
		\begin{align*}
			\delta_{m, \left\lceil k/2 \right\rceil } 	  	& \geq f_{0,k}  			\\ 
			\delta_{m, \left\lceil (n+k)/2 \right\rceil } 	& \geq -f_{n-k-1,1}  		\\ 
			\epsilon_{m, a-b} 								& \geq f_{a,b} - f_{a,b+1} 	\\
			\epsilon_{m, a-b} 								& \geq f_{a,b} - f_{a-1,b}.
		\end{align*}
		If $m = \left\lceil k/2 \right\rceil$, then this is precisely the promised Gelfand-Tsetlin polytope. Otherwise we need a translation by an integer vector depending on $m$.
		For $m = \left\lceil (n+k)/2 \right\rceil$ the translation
		$$f_{a,b} \mapsto f_{a,b}+1$$
		works.
		Now if $m$ is neither $\left\lceil k/2 \right\rceil$ nor $\left\lceil (n+k)/2 \right\rceil$, then $\epsilon_{m, d}$ evaluates to $1$ for exactly one $d \in \{1-k, \ldots, n-k-2\}$. Then translation 
		$$ f_{a,b} \mapsto f_{a,b} + \delta_{a-b<d}, $$
		where $\delta_{a-b<d} = 1$ iff $a-b < d$, does the trick.
	\end{proof}
	
	Up to some minor changes, the above proof also works for any graph in the orbit of the checkboard graphs as well as the dual checkboard graphs. This mostly boils down to adjusting the map $\tilde F$ and the translations at the end.
	
	\begin{theorem}
		Let $G$ be any graph in the $D_n$-orbit of either $G^\text{ch}_{k,n}$ or $\widehat G^\text{ch}_{k,n}$. Then $\Delta_G$ is unimodular equivalent to the Gelfand-Tsetlin polytope $\mathcal{GT}_{k,n}^1$.
	\end{theorem}

	\bibliographystyle{apalike}
	\bibliography{NewtonOkounkovBodiesObtainedFromCertainOrbitsOfPlabicGraphs}
	
\end{document}